\documentclass[preprint,11pt]{elsarticle}

\usepackage{amssymb}
\usepackage{amsthm}

\usepackage{multicol}
\usepackage{mathrsfs}
\usepackage{amsfonts}
\usepackage{CJK,fancybox,color,graphicx,amsmath,amsthm,amssymb,enumerate}

\newtheorem{thm}{Theorem}[section]

\newtheorem{prop}[]{Proposition}[section]
\theoremstyle{definition}

\theoremstyle{remark}
\newtheorem{rem}[]{Remark}[section]
\numberwithin{equation}{section}

\newcommand{\D}{\displaystyle}
\newcommand{\DF}[2]{\frac{\D#1}{\D#2}}

\begin{document}
\begin{frontmatter}

\title{Singular solutions to the Riemann problem for the pressureless Euler equations with discontinuous source term}


\author{Qingling Zhang}

\address{School of Mathematics and Computer Sciences, Jianghan University, Wuhan 430056, PR China }
\ead{zhangqingling2002@163.com}
\begin{abstract}
In this paper, the Riemann problem for the pressureless Euler equations with a discontinuous source term is considered.
The delta shock wave solution is obtained by combining the generalized Rankine-Hugoniot conditions together with the method of
characteristics for different situations, which reflects the impact of the source term on the delta shock front. Moreover, during the construction process
of the Riemann solution, some
interesting phenomena are also observed, such as the disappearance of the delta shock wave and the occurrence of the vacuum state, etc.

\end{abstract}


\end{frontmatter}

\section{Introduction}
It is well know that the singular discontinuities ( delta shock waves) may develop for a hyperbolic conservation laws, which
may relult from the initial data or the degeneracy and coincidence of the characteristics.
The propagation of such delta shock waves or singular discontinuities may be affected by source terms for nonlinear systems, such as
 the pressureless Euler system with friction or damping \cite{Shen1} and the (generalized) Chaplygin gas with friction term \cite{Shen2,Sun2}.
In this paper, we are concerned with the pressureless Euler system with discontinuous source term
in the following form:
\begin{equation}\label{1.1}
\left\{\begin{array}{ll}
\rho_t+(\rho u)_x=0,\\
u_t+(\frac{u^2}{2})_x=H(x-s(t))f(x,t,u)+H(s(t)-x)g(x,t,u),
\end{array}\right.
\end{equation}
 where $H$ is the standard Heaviside function, $f(x,t,u)$ and $g(x,t,u)$ is a given continuous
 function with respect to $x$ and $t$, and $\rho>0$ and $u$ denote the density and the velocity, respectively.
From the second equation of (\ref{1.1}), $u$ is assumed to be a piecewise smooth function with a
single jump discontinuity at the curve $x=s(t)$ in the $(x,t)$ plane.

There are also quite a number of
physical phenomena that can be described by hyperbolic conservation laws with singular source terms, such as
the shallow-water flow of gravity currents \cite{Montgomery-Moodie1}, atmospheric cold fronts \cite{Montgomery-Moodie2} and radiative hydrodynamics \cite{Mihalas-Mihalas}. Moreover,
the regularization technique has been developed in \cite{Suarez-Jacobs-Don} to deal with hyperbolic conservation laws with time-dependent singular Dirac delta source terms. In the future, we will also consider the pressureless Euler system with singular source term.

The main purpose of this paper is to consider the impact of the inhomogeneous source term $f(x,t,u)$ and $g(x,t,u)$ on
the location of delta shock front. In order to give a explicit expression of the delta shock wave curve
and to display the effect of the source term, we will take $f(x,t,u)=0$, $g(x,t,u)=1$ or $-u$ and $f(x,t,u)=1$, $g(x,t,u)=-u$ as typical examples
to study the system (\ref{1.1}) with Riemann initial data
\begin{equation}\label{1.2}
(\rho,u)(x,0)=\left\{\begin{array}{ll}
(\rho_-,u_-),\ \ x<0,\\
(\rho_ +,u_+),\ \ x>0.
\end{array}\right.
\end{equation}
With the above discontinuous source terms, the delta shock wave solutions can display some interesting behaviors and the vacuum state occurs in some situations,
which may provide some insights into more general source term situations.

Specifically, for $f(x,t,u)=0$, $g(x,t,u)=1$ or $-u$, (\ref{1.1}) is reduce to the following form:
\begin{equation}\label{1.3}
\left\{\begin{array}{ll}
\rho_t+(\rho u)_x=0,\\
u_t+(\frac{u^2}{2})_x=H(s(t)-x)g(x,t,u),
\end{array}\right.
\end{equation}
For smooth solutions, (\ref{1.3}) is equivalen to the following system:
\begin{equation}\label{1.4}
\left\{\begin{array}{ll}
\rho_t+(\rho u)_x=0,\\
(\rho u)_t+(\rho u^2)_x=H(s(t)-x)\rho g(x,t,u),
\end{array}\right.
\end{equation}
which is just the generalization of the pressureless Euler system with friction or damping \cite{Ha-Huang-Wang,Jin,Shen1}.
Recently, we discussed the stability of the Rimann solutions to the the pressureless Euler system with Coulomb-like friction by flux approximation in \cite{Zhang}.
For other cases such as $g(x,t,u)=x$ corresponding to \cite{Ding-Huang}, we will consider it in the future.

Let us assume $u_{-}>u_{+}$, then a singular discontinuous solution $(\rho,u)$ is desired to be found
in the domain $-\infty<x<\infty$ and $t>0$. For the homogeneous pressureless Euler system, there are already many results, see \cite{Huang-Wang,Shen-Sun1,Sheng-Zhang,Wang-Ding,Wang-Huang-Ding}
For the homogeneous case $f(x,t,u)=g(x,t,u)=0$ with $u_{-}>u_{+}$, it is well known that the Riemann problem  owns
a delta shock solution, which should satisfy the following generalized Rankine-Hugoniot condition \cite{Shen-Sun1}:
\begin{equation}\label{1.5}
\left\{\begin{array}{ll}
\DF{ds(t)}{dt}=\sigma(t)=u_{\delta},\\[4pt]
\DF{dw(t)}{dt}=\sigma(t)[\rho]-[\rho u],\\[4pt]
\sigma(t)[u]=[\frac{u^{2}}{2}],
\end{array}\right.
\end{equation}
whose Riemann solution can be expressed as
\begin{equation}\label{1.6}
(\rho,u)(x,t)=\left\{\begin{array}{ll}
(\rho_-,u_-),\ \ x<x(t),\\
(w(t),u_{\delta}(t)),\ \ x=x(t)\\
(\rho_ +,u_+),\ \ x>(t).
\end{array}\right.
\end{equation}
where $x=s(t)=\frac{u_{-}+u_{+}}{2}t$,
$w(t)=\frac{1}{2}(\rho_-+\rho_+)(u_{-}-u_{+})$ and
$\frac{dx(t)}{dt}=\sigma(t)=\frac{u_{-}+u_{+}}{2}$ are respectively
denote the location, weight and propagation speed of the delta shock
wave, and $u_{\delta}(t)$ indicates the assignment of $u$ on this
delta shock wave, and $[\rho]=\rho(x(t)+0,t)-\rho(x(t)-0,t)$ denotes
the jump of the function $\rho$ across the delta shock wave.

For the inhomogeneous situations, it is worhtwhile to notice that the characteristic curves do not keep  straight and
$u$ along each of the characteristic curves is not a constant again on one side or both sides of the delta shock front.
Furthermore, the source term will act on the delta shock front such that it will also bend and no longer be a straight line again.
Thus, the characteristic curves and the values of $u$ along these characteristic curves should be solved by applying the method of characteristics, and then the generalized Rankine-Hugoniot conditions should be proposed to determine the location of the delta shock wave front. For the Riemann problem for
hyperbolic conservation laws, the readers can refer to the standard textbooks such as \cite{Chang-Hsiao,Smoller}.
Recently, there are many works concentrate on how to use the method of characteristic to solve the Riemann problem for scalar conservation law with discontinuous coefficient or with source term, such as \cite{Chang-Chou-Hong-Lin, Fang-Tang-Wang,Sun1}. However, there are very few related works for conservation laws, such as \cite{Tang-Wang}.
In this paper and in the future work, we will focus on Riemann problem for conservation laws with discontinuous coefficients or with source terms.

It is clear to see that (\ref{1.1}) depends on the existence of the delta shock wave. Thus, if the delta shock wave
does not exist or the entropy condition cannot hold, then we shall cease to construct the delta shock wave solution.
In fact, a delta shock wave is generated for the Riemann problem (\ref{1.1}) and (\ref{1.2}) provided that $u_{-}>u_{+}$.
Some interesting phenomena can be observed about the structural behavior of the delta shock wave in the Riemann solution to (\ref{1.1}) and (\ref{1.2}).

Obviouly, the second equation of (\ref{1.1}) is a special type of the quasilinear hyperbolic equations
\begin{equation}\label{1.7}
u_t+F(x,t,u)_x=G(x,t,u).
\end{equation}
The weak solution to the Cauchy problem to the scalar situation for (\ref{1.7}) was first studied in \cite{Kruzkov} and \cite{Volpert}. For the general form $F(x,t,u)$ and $G(x,t,u)$, the situations are more complicated and still far from
a complete description of Riemann solutions. Therefore, some special choices of $F(x,t,u)$ and $G(x,t,u)$ are taken to
study the possible structures and asymptotic profiles of solutions. The particular Riemann
problem (\ref{1.5}) and (\ref{1.2}) has been extensively studied, such as in \cite{Sinestrari} for $F(x,t,u)=F(u)$ and $G(x,t,u)=G(u)$,
in \cite{Greenberg-Leroux,Hong-Temple,Karlsen-Risebro-Towers} for $F(x,t,u)=F(u)$ and $G(x,t,u)=a'(x)u$, and
in \cite{Greenberg-Leroux-Baraille,Karlsen-Mishra-Risebro} for $F(x,t,u)=F(u)$ and $G(x,t,u)=a'(x)$, where $a(x)$ is usually assumed to be discontinuous at $x=0$
and its derivative is also a Dirac delta measure such that a standing wave discontinuity is developed on the line $x=0$.
For more cases about the hyperbolic conservation laws with source terms, the author can see in \cite{Lien,Liu} for $F(x,t,u)=F(u)$ and $G(x,t,u)=c(x-\alpha t)g(u)$
where $\alpha$ is the speed of the source, in \cite{Diehl} for $F(x,t,u)=F(x,u)$ and $G(x,t,u)=s(t)\delta(x)$ where the source models the inlet,etc. The author can also refer to \cite{E-Khanin-Mazel-Sinai,Saichev-Woyczynski} for the related results about the Burgers equation with source terms.

Recently, Fang et al. \cite{Fang-Tang-Wang} have investigated the Riemann problem for the Burgers equation with a discontinuous source term as follows:
\begin{equation}\label{1.8}
u_t+(\frac{u^2}{2})_x=g(x,t),
\end{equation}
where $g(x,t)=g^-H(-x)+g^+H(x)$ and $g^-,g^+$
are two constants.
The Riemann problem (\ref{1.8}) and (\ref{1.2}) have been constructed completely in \cite{Fang-Tang-Wang} where some interesting phenonmena were observed there.
Moreover, they have also studied the behavior of the shock waves in a $2\times2$
balance law with discontinuous source terms.
In this paper, similar to \cite{Zhang-Shen}, we assume that the
discontinuity of the source term alwalys accompanies with the location of the delta shock front, so we only need to concern the delta
shock wave solution to the Riemann problem (\ref{1.1}) and (\ref{1.2}). Moreover, we find that this assumption in our paper simplifies our study greatly
which enable us to explicitly construct the delta shock wave solution to Riemann problem (\ref{1.1}) and (\ref{1.2})
 and furthermore clearly display the effect of the source term on the delta shock wave.

The paper is organized as follows. In Section 2, the generalized Rankine-Hugoniot condition for (\ref{1.1}) is derived in order to be selfcontained.
In Section 3, the Riemann problem for (\ref{1.1}) and (\ref{1.2}) is considered and the delta shock wave solution is constructed when
$u_{-}>u_{+}$, $f(x,t,u)=0$ and $g(x,t,u)=1$ or $-u$.  In Section 4, the delta shock wave solution to the Riemann problem (\ref{1.1}) and (\ref{1.2}) is constructed when
$u_{-}>u_{+}$, $f(x,t,u)=1$ and $g(x,t,u)=-u$. Finally, the discussions are carried our and the conclusions are drawn in Section 5.

\section{The generalized Rankine-Hugoniot jump conditions}
For a nonlinear system of hyperbolic conservation laws with a source term which is discontinuous at the location of
the shock front, Montgomery and Moodie \cite{Montgomery-Moodie1} proposed a generalized Rankine-Hugoniot jump condition to determine the shock front.
For the case when the source term is continuous, Shen \cite{Shen1}give the generalized Rankine-Hugoniot jump condition for the delta shock wave front.
In this paper, we only restrict ourselves to the particular form (\ref{1.1}) for the convenience of study.

\begin{prop}\label{prop:2.1}
Assume that
\begin{equation}\label{2.1}
(\rho,u)(x,t)=\left\{\begin{array}{ll}
(\rho_l,u_l),\ \ x<s(t),\\
(w(t),u_{\delta}(t)),\ \ x=s(t)\\
(\rho_r,u_r),\ \ x>s(t),
\end{array}\right.
\end{equation}
is a delta shock wave solution of (\ref{1.1}) which has a single jump discontinuity at the
position of the delta shock front $x=s(t)$. Then, at $x=s(t)$, the jump conditions (\ref{1.5}) still hold
and the following $\delta$-entropy condition is satisfied:
\begin{equation}\label{2.2}
u(s(t)+0)<\frac{ds}{dt}<u(s(t)-0).
\end{equation}
\end{prop}

\textbf{Proof}  We will only prove the third identity of
(\ref{1.5}). For the proof of the others, one can refer to \cite{Shen-Sun1}. Let us
assume that $\Gamma$ is the delta shock wave curve. Then, we use
$\Omega_{-}:{(x,t)x<s(t)}$ to denote the left-hand side region of
$\Gamma$ and $\Omega_{+}:\{(x,t)\mid x>s(t)\}$ to denote the
right-hand side region of $\Gamma$. It is clear to see from
(\ref{1.1}) that we have
\begin{equation}\label{2.3}
u_t+(\frac{u^2}{2})_x=g(x,t,u),\                                                                                                                                                                                                                                                                                                     \ in \Omega_-,
\end{equation}
\begin{equation}\label{2.4}
u_t+(\frac{u^2}{2})_x=f(x,t,u),\ \ in \Omega_+.
\end{equation}

Let $\varphi(t)\in C_{c}^{\infty}(R_{+}^{2})$ be a test function, such that the following equality
\begin{equation}\label{2.5}
\langle u_t+(\frac{u^2}{2})_x,\varphi\rangle=\langle H(x-s(t))f(x,t,u)+H(s(t)-x)g(x,t,u),\varphi\rangle
\end{equation}
should hold in the sense of distributions. It is clear to see that the left-hand side of (\ref{2.5})
can be calculated by
\begin{eqnarray}
\langle u_t+(\frac{u^2}{2})_x,\varphi\rangle&=&-\langle u,\varphi_t\rangle-\langle\frac{u^2}{2},\varphi_x\rangle\nonumber\\
&=&-\int\int_{\Omega_{-}}\Big(u\varphi_{t}+\frac{u^2}{2}\varphi_{x}\Big)dx
dt-\int\int_{\Omega_{+}}\Big(u\varphi_{t}+\frac{u^2}{2}\varphi_{x}\Big)dx
dt\nonumber\\
&=&-\int\int_{\Omega_{-}}\Big((u\varphi)_{t}+(\frac{u^2}{2}\varphi)_{x}\Big)dx
dt+\int\int_{\Omega_{-}}g\varphi dx
dt\nonumber\\&&-\int\int_{\Omega_{+}}\Big((u\varphi)_{t}+(\frac{u^2}{2}\varphi)_{x}\Big)dx
dt+\int\int_{\Omega_{+}}f\varphi dx
dt\nonumber\\
&=&\int_{\partial\Omega_{-}}\frac{u^2}{2}\varphi dt-u\psi
dx+\int\int_{\Omega_{-}}g\varphi dx
dt\nonumber\\&&\int_{\partial\Omega_{+}}\frac{u^2}{2}\varphi
dt-u\psi dx+\int\int_{\Omega_{+}}f\varphi dx
dt\nonumber\\
&=&\int_{0}^{\infty}(\sigma[u]-[\frac{u^2}{2}])\varphi
dt+\int\int_{\Omega_{-}}g\varphi dx dt+\int\int_{\Omega_{+}}f\varphi
dx dt\label{2.6},
\end{eqnarray}
in which (\ref{2.3}) and (\ref{2.4}) have been used.

On the other hand, the right-hand side of (\ref{2.5}) can be
calculated by
\begin{equation}\label{2.7}
\langle
H(x-s(t))f(x,t,u)+H(s(t)-x)g(x,t,u),\varphi\rangle=\int\int_{\Omega_{-}}g\varphi
dx dt+\int\int_{\Omega_{+}}f\varphi dx dt.
\end{equation}

 For $\varphi$ is arbitrary, then the third identity of
(\ref{1.5}) can be obtained by combining (\ref{2.6}) with
(\ref{2.7}) together.

\section{The situation for $f(x,t,u)=0$ or $g(x,t,u)=0$}
Without loss of generality, we assume that $f(x,t,u)=0$ and
$g(x,t,u)\neq0$. Specially, for simplicity and practice meaning in
\cite{Ha-Huang-Wang,Jin,Shen1}, we take $g(x,t,u)=1$ or $-u$ in the following two subsections.

\subsection{The situation for $g(x,t,u)=1$}
In this subsection, we consider the Riemann problem (\ref{1.1}) and
(\ref{1.2}) for the situation $f(x,t,u)=0$ and $g(x,t,u)=1$. When
$u_{-}>u_{+}$, it is clear to see that the Riemann problem
(\ref{1.1}) and (\ref{1.2}) has a delta shock solution with a single
jump discontinuity across the path $x=s(t)$, which crosses the
$x-$axis at $s(0)=0$ and satisfies the generalized Rankine-Hongniot
condition (\ref{1.5}). So the the Riemann problem (\ref{1.1}) and
(\ref{1.2}) can be found by considering the prolem on either side of
the delta shock using the method of characteristics and then
implementing the generalized Rankine-Hongniot condition (\ref{1.5})
to determine the position of the delta shock front.

\begin{thm}\label{thm:3.1}
For the situation $f(x,t,u)=0$ and $g(x,t,u)=1$, if $u_{-}>u_{+}$,
then the Riemann solution to (\ref{1.1}) and (\ref{1.2}) is a delta
shock solution which can be express as
\begin{equation}\label{3.1}
(\rho,u)(x,t)=\left\{\begin{array}{ll}
(\rho_-,u_-+t),\ \ x<s(t),\\
(w(t),u_{\delta}(t)),\ \ x=s(t)\\
(\rho_+,u_+),\ \ x>s(t),
\end{array}\right.
\end{equation}
with
\begin{equation}\label{3.2}
 s(t)=\frac{1}{2}t^{2}+\frac{1}{2}(u_{-}+u_{+})t,
\end{equation}
\begin{equation}\label{3.3}
 w(t)=\frac{1}{2}(\rho_{-}+\rho_{+})\{\frac{1}{2}t^{2}+(u_{-}-u_{+})t\},
\end{equation}
\begin{equation}\label{3.4}
 u_{\delta}(t)=\frac{1}{2}(u_{-}+u_{+}+t).
\end{equation}
Furthermore, the delta shock exists for all the time by taking into account the entropy condition.

\end{thm}

\textbf{Proof}  \hspace{0.5cm}If $x>s(t)$, then we have $H(s(t)-x)=0$ and the
initial value problem (\ref{1.1}) and (\ref{1.2}) becomes
\begin{equation}\label{3.5}
\left\{\begin{array}{ll}
\rho_t+(\rho u)_x=0,\\
u_t+uu_x=0,
\end{array}\right.
\end{equation}
with initial condition
\begin{equation}\label{3.6}
(\rho,u)(x_{0},0)=(\rho_+,u_+),where \ \ x_0>s(0)=0.
\end{equation}

It is obvious to see that the problem (\ref{3.5}) and (\ref{3.6})
has a trivial solution as
\begin{equation}\label{3.7}
(\rho,u)(x,t)=(\rho_+,u_+), for all \ \ x>s(t), t\geq0.
\end{equation}
Thus, if $x_0>0$, then the characteristic curve starting from the
initial point $(x_0,0)$ on the $x$-axis can be calculated by
\begin{equation}\label{3.8}
\frac{dx}{dt}=u, with \ \ x\mid_{t=0}=x_0,
\end{equation}
which implied that
\begin{equation}\label{3.9}
x=u_+t+x_0.
\end{equation}

On the other hand, if $x<s(t)$, then we have $H(s(t)-x)=1$ and the
initial value problem (\ref{1.1}) and (\ref{1.2}) becomes
\begin{equation}\label{3.10}
\left\{\begin{array}{ll}
\rho_t+(\rho u)_x=0,\\
u_t+uu_x=1,
\end{array}\right.
\end{equation}
with initial condition
\begin{equation}\label{3.11}
(\rho,u)(x_{0},0)=(\rho_-,u_-),where \ \ x_0<s(0)=0.
\end{equation}

By the method of characteristics, it is obvious to see that the
problem (\ref{3.10}) and (\ref{3.11}) has a trivial solution as
\begin{equation}\label{3.12}
(\rho,u)(x,t)=(\rho_-,u_-+t), for all \ \ x<s(t), t\geq0.
\end{equation}
Thus, if $x_0<0$, then the characteristic curve starting from the
initial point $(x_0,0)$ on the $x$-axis can be calculated by
\begin{equation}\label{3.13}
\frac{dx}{dt}=u_-+t, with \ \ x\mid_{t=0}=x_0,
\end{equation}
which implied that
\begin{equation}\label{3.14}
x=\frac{1}{2}t^{2}+u_-t+x_0.
\end{equation}

To connect the solutions of (\ref{3.7}) and (\ref{3.12}) as a delta
shock solution to the Riemann problem to (\ref{1.1}) and
(\ref{1.2}), the generalized Rankine-Hongniot condition (\ref{1.5})
should be imposed at the position of the delta shock front $x=s(t)$
as follows:
\begin{equation}\label{3.15}
\frac{ds}{dt}(u_+-(u_-+t))=\frac{1}{2}({u_+}^2-(u_-+t)^2),
\end{equation}
namely, we have
\begin{equation}\label{3.16}
u_{\delta}(t)=\sigma(t)=\frac{ds}{dt}=\frac{1}{2}({u_-}+u_++t).
\end{equation}
With $s(0)=0$ in mind, we can get the expression of the delta shock
front (\ref{3.2}).

Moreover, the weight of the delta shock can be get from the second
equality of (\ref{1.5}) with
\begin{eqnarray}
\frac{d
w(t)}{dt}&=&\sigma(t)(\rho_+-\rho_--(\rho_+u_+-\rho_-(u_-+t)))\nonumber\\
&=&\frac{1}{2}({u_-}+u_++t)(\rho_+-\rho_--(\rho_+u_+-\rho_-(u_-+t)))\nonumber\\
&=&\frac{1}{2}(\rho_++\rho_-)(u_--u_++t)\label{3.17},
\end{eqnarray}
from which we get (\ref{3.3}).

Since the $\delta$-entropy condition
\begin{equation}\label{3.18}
u_{+}<\frac{ds}{dt}<u_{-}+t, for t\geq0,
\end{equation}
is always satisfied, so the the delta shock always exists. The proof is completed.

\bigbreak
Now we are in the position to consider the path of the delta shock wave for
the Riemann problem (\ref{1.1}) and (\ref{1.2}). For the nomogeneous situation, the
he path of the delta shock wave is a straight line with the slope of $\frac{u_-+u_+}{2}$.
For the innomogeneous situation with $g(x,t,u)=1$, the source term will make the delta shock front bend.
From (\ref{3.16}) we have
\begin{equation}\label{3.19}
\frac{ds}{dt}\mid_{t=0}=\frac{1}{2}({u_-}+u_+),
\end{equation}
which implies that the speed of the delta shock wave for innomogeneous situation with $g(x,t,u)=1$ is identical with that of
the nomogeneous situation in the beginning.

In order to describe the effect of the forcing term in detail, the discussion should be divided into the following two cases.
 \bigbreak
\noindent Case 1 \hspace{0.5cm} If ${u_-}+u_+\geq 0$, taking into account to ${u_-}>u_+$, it is easy to get $u_->0$. Moreove, we can get
$\frac{ds}{dt}>0$ and $\frac{d^{2}s}{dt^{2}}=\frac{1}{2}>0$ for any $t>0$ from (\ref{3.16}). In other words, $s(t)$
is always convex and increases along with $t$
such that the delta shock wave should always move forward.
Let us draw Fig.1(a) and Fig.1(b) for the situations ${u_-}>u_+>0$ and ${u_-}>0\geq u_+$ respectively.

\bigbreak
\noindent Case 2 \hspace{0.5cm} If ${u_-}+u_+<0$, taking into account to ${u_-}>u_+$, it is easy to get $u_+<0$,
and from (\ref{3.19}) the delta shock wave has a negative speed in the beginning.
Since $\frac{d^{2}s}{dt^{2}}=\frac{1}{2}>0$ for any $t>0$ from (\ref{3.16}), $s(t)$
is always convex, which implies that the speed of the delta shock wave speeds up. It is clear that there exists a time $t_1=-\frac{1}{2}({u_-}+u_+)$
, such that $\frac{ds}{dt}\mid_{t=t_1}=0$. In other words, the delta shock moves backword for $0<t<t_1$
and moves forward for $t>t_1$. It follows from (\ref{3.2}) that there also extists a time $t_2=-({u_-}+u_+)$
such that $s(t_{2})=0$, which means that the delta shock front intersects with the $t$-axis at the time $t_2$.
Let us draw Fig.1(c) and Fig.1(d) for the situations $u_-\geq0>u_+$ and $0>u_{-}>u_+$ respectively, where $t=t_1^\ast$ is the the symmetry axis of the characteristic
curves on the right hand side of the delta shock waves.

\begin{multicols}{2}
\unitlength 1mm 
\linethickness{0.4pt}
\ifx\plotpoint\undefined\newsavebox{\plotpoint}\fi 
\begin{picture}(87.25,48.5)(15,0)
\put(86.25,7.5){\vector(1,0){.07}}
\put(10.25,7.5){\line(1,0){75}}
\put(46,48.5){\vector(0,1){.07}}
\put(46,7.5){\line(0,1){40.5}}
\qbezier(46,7.5)(52.75,32.75)(75.5,44.75)
\qbezier(53.5,25.25)(40.13,20)(31.25,7.5)
\qbezier(62.25,35.75)(34.13,22.25)(23.5,7.5)
\qbezier(70.25,42)(30.63,29.75)(16.5,7.5)
\put(43.75,48){$t$}
\put(87.25,5.25){$x$}
\put(46.25,4.75){0}
\put(36.75,2){(a) $u_->u_+>0$}
\put(70.25,45.75){$\delta S$}

\multiput(53.5,25.25)(-.03348214,-.15848214){112}{\line(0,-1){.15848214}}
\put(57.25,7.5){\line(0,1){0}}
\multiput(62.25,35.5)(-.03370787,-.15589888){180}{\line(0,-1){.15589888}}
\put(68.25,7.75){\line(0,1){0}}
\multiput(70.75,42)(-.030073469,-.139795918){248}{\line(0,-1){.139795918}}
\put(79,7.75){\line(0,1){0}}
\end{picture}

\unitlength 1mm 
\linethickness{0.4pt}
\ifx\plotpoint\undefined\newsavebox{\plotpoint}\fi 
\begin{picture}(87.25,48.5)(10,0)
\put(86.25,7.5){\vector(1,0){.07}}
\put(10.25,7.5){\line(1,0){75}}
\put(46,48.5){\vector(0,1){.07}}
\put(46,8){\line(0,1){40.5}}
\qbezier(46,7.5)(52.75,32.75)(75.5,44.75)
\qbezier(53.5,25.25)(40.13,20)(31.25,7.5)
\qbezier(62.25,35.75)(34.13,22.25)(23.5,7.5)
\qbezier(70.25,42)(30.63,29.75)(16.5,7.5)
\put(43.75,48){$t$}
\put(87.25,5.25){$x$}
\put(46.25,4.75){0}
\put(36.75,2){(b) $u_->0\geq u_+$}
\put(70.25,45.75){$\delta S$}
\multiput(53.5,25.25)(.03348214,-.15848214){112}{\line(0,-1){.15848214}}
\put(57.25,7.5){\line(0,1){0}}
\multiput(62.25,35.5)(.03370787,-.15589888){180}{\line(0,-1){.15589888}}
\put(68.25,7.75){\line(0,1){0}}
\multiput(70.75,42)(.033673469,-.139795918){248}{\line(0,-1){.139795918}}
\put(79,7.75){\line(0,1){0}}
\end{picture}
\end{multicols}


\begin{multicols}{2}
\unitlength 0.9mm 
\linethickness{0.4pt}
\ifx\plotpoint\undefined\newsavebox{\plotpoint}\fi 
\begin{picture}(98.5,55.5)(15,0)
\put(43,55){$t$}
\put(97.25,4.25){$x$}
\put(46.25,4){0}
\put(32,10.5){t=$t_1$}
\put(39,11.5){$\bullet$}
\put(36.75,1){(c) $u_-\geq 0>u_+$}
\put(70.25,45.75){$\delta S$}
\put(57.25,7.5){\line(0,1){0}}
\put(68.25,7.75){\line(0,1){0}}
\put(79,7.55){\line(0,1){0}}
\qbezier(45.75,7.5)(25.13,14.75)(79,48)
\multiput(44,21.25)(.03372093,-.03372093){422}{\line(0,-1){.03372093}}
\put(58.5,7.5){\line(0,1){0}}
\multiput(54,30.5)(.0337243402,-.0337243402){682}{\line(0,-1){.0337243402}}
\put(77,7.5){\line(0,1){0}}
\put(77,7.5){\line(0,1){0}}
\multiput(63.5,37.75)(.0337142857,-.0337142857){890}{\line(0,-1){.0337142857}}
\put(93,8.25){\line(0,1){0}}
\qbezier(63,38)(23.63,26.63)(11.75,7.5)
\put(45.5,55.5){\vector(0,1){.07}}
\put(45.5,7.5){\line(0,1){47.75}}
\qbezier(44,21)(37.25,19.5)(27.5,7.5)
\qbezier(54,30.75)(29.75,21.75)(19.5,7.5)
\put(98.5,7){\vector(1,0){.07}}
\put(8,7.5){\line(1,0){90.5}}
\end{picture}
\unitlength 0.9mm 
\linethickness{0.4pt}
\ifx\plotpoint\undefined\newsavebox{\plotpoint}\fi 
\begin{picture}(98.5,55.5)(5,0)

\put(44.4,22){$\bullet$}
\put(48,22){$t_2$}
\put(43,55){$t$}
\put(97.25,4.25){$x$}
\put(46.25,4){0}
\put(36.75,1){(d) $0>u_-> u_+$}
\put(70.25,45.75){$\delta S$}
\put(57.25,7.5){\line(0,1){0}}
\put(68.25,7.75){\line(0,1){0}}
\put(79,7.75){\line(0,1){0}}
\qbezier(45.75,7.75)(25.13,14.75)(79,48)
\multiput(44,21.25)(.03372093,-.03372093){411}{\line(0,-1){.03372093}}
\put(58.5,6.75){\line(0,1){0}}
\multiput(54,30.5)(.0337243402,-.0337243402){680}{\line(0,-1){.0337243402}}
\put(77,7.75){\line(0,1){0}}
\put(77,6.5){\line(0,1){0}}
\multiput(63.5,37.75)(.0337142857,-.0337142857){879}{\line(0,-1){.0337142857}}
\put(93,7.75){\line(0,1){0}}
\put(45.5,55.5){\vector(0,1){.07}}
\put(45.5,7.75){\line(0,1){47.75}}
\put(45.5,55.5){\vector(0,1){.07}}
\put(45.5,55.5){\line(0,1){0}}
\put(98.5,7.75){\vector(1,0){.07}}
\put(8,7.75){\line(1,0){90.5}}

\bezier{50}(8,9.5)(18,9.5)(38.5,9.5)
\put(4,10.5){t=$t_1^\ast$}

\put(32,11){t=$t_1$}
\put(39,11.5){$\bullet$}
\put(11.5,8.5){$\bullet$}
\put(20,8.5){$\bullet$}
\put(28,8.5){$\bullet$}
\qbezier(43.75,21.5)(20.13,9.75)(34,7.75)
\qbezier(53.5,31)(10.75,9.63)(24,7.75)
\qbezier(63.75,38)(.38,10.88)(15.5,7.75)
\put(12.25,-3){\makebox(0,0)[cc]
{Fig.1 The delta shock wave solution to (1.1) and (1.2) when $f(x,t,u)=0$ and $g(x,t,u)=1$ }}

\end{picture}
\end{multicols}

\bigbreak
\begin{rem}\label{rem:3.1}   Similarly, we can consider the situation $g(x,t,u)=-1$.
\end{rem}

\subsection{The situation for $g(x,t,u)=-u$}

In this subsection, we consider the Riemann problem (\ref{1.1}) and
(\ref{1.2}) for the situation $f(x,t,u)=0$ and $g(x,t,u)=-u$. Similar to subsection 3.1,
we have the following theorem to depict the delta shock wave solution to the Riemann problem (\ref{1.1}) and
(\ref{1.2}) when $u_{-}>u_{+}$.

\begin{thm}\label{thm:3.2}
For the situation $f(x,t,u)=0$ and $g(x,t,u)=-u$, if $u_{-}>u_{+}$,
then the Riemann solution to (\ref{1.1}) and (\ref{1.2}) is a delta
shock solution which can be express as
\begin{equation}\label{3.20}
(\rho,u)(x,t)=\left\{\begin{array}{ll}
(\rho_-,u_-e^{-t}),\ \ x<s(t),\\
(w(t),u_{\delta}(t)),\ \ x=s(t)\\
(\rho_+,u_+),\ \ x>s(t),
\end{array}\right.
\end{equation}
with
\begin{equation}\label{3.21}
 s(t)=\frac{1}{2}(u_{-}(1-e^{-t})+u_{+}t),
\end{equation}
\begin{equation}\label{3.22}
 w(t)=\frac{1}{2}(\rho_{-}+\rho_{+})(u_{-}(1-e^{-t})-u_{+}t),
\end{equation}
\begin{equation}\label{3.23}
 u_{\delta}(t)=\frac{1}{2}(u_{-}e^{-t}+u_{+}).
\end{equation}
Furthermore, the delta shock exists for all the time by taking into account the entropy condition.

\end{thm}

\textbf{Proof}  \hspace{0.5cm}If $x>s(t)$, then we have $H(s(t)-x)=0$. Similar to subsection 3.1, we have
\begin{equation}\label{3.24}
(\rho,u)(x,t)=(\rho_+,u_+), for all \ \ x>s(t), t\geq0.
\end{equation}
And if $x_0>0$, then the characteristic curve starting from the
initial point $(x_0,0)$ on the $x$-axis is
\begin{equation}\label{3.25}
x=u_+t+x_0.
\end{equation}

On the other hand, if $x<s(t)$, then we have $H(s(t)-x)=1$ and the
initial value problem (\ref{1.1}) and (\ref{1.2}) becomes
\begin{equation}\label{3.26}
\left\{\begin{array}{ll}
\rho_t+(\rho u)_x=0,\\
u_t+uu_x=-u,
\end{array}\right.
\end{equation}
with initial condition
\begin{equation}\label{3.27}
(\rho,u)(x_{0},0)=(\rho_-,u_-),where \ \ x_0<s(0)=0.
\end{equation}

By the method of characteristics, it is obvious to see that the
problem (\ref{3.26}) and (\ref{3.27}) has a trivial solution as
\begin{equation}\label{3.28}
(\rho,u)(x,t)=(\rho_-,u_-e^{-t}), for all \ \ x<s(t), t\geq0.
\end{equation}
Similarly, if $x_0<0$, then the characteristic curve starting from the
initial point $(x_0,0)$ on the $x$-axis can be calculated by
\begin{equation}\label{3.29}
\frac{dx}{dt}=u_-e^{-t}, with \ \ x\mid_{t=0}=x_0,
\end{equation}
which implied that
\begin{equation}\label{3.30}
x=u_-(1-e^{-t})+x_0.
\end{equation}

It is clear to see that for some situation, after some time $t$ the delta shock wave will disappear for the reason that
the $\delta$-entropy condition cannot be satisfied. However, for a sufficiently small time $t$, the delta shock wave may exist and should satisfy
the generalized Rankine-Hongniot condition (\ref{1.5}), so we have
\begin{equation}\label{3.31}
\frac{ds}{dt}(u_+-u_-e^{-t} )=\frac{1}{2}({u_+}^2-(u_-e^{-t})^2),
\end{equation}
namely, we have
\begin{equation}\label{3.32}
u_{\delta}(t)=\sigma(t)=\frac{ds}{dt}=\frac{1}{2}({u_-}e^{-t}+u_+).
\end{equation}
So (\ref{3.23})is obtained. With $s(0)=0$ in mind, we can get the expression of the delta shock
front (\ref{3.21}).

Moreover, the weight of the delta shock can be get from the second
equality of (\ref{1.5}) with
\begin{eqnarray}
\frac{d
w(t)}{dt}&=&\sigma(t)(\rho_+-\rho_-)-(\rho_+u_+-\rho_-u_-e^{-t})\nonumber\\
&=&\frac{1}{2}({u_-}e^{-t}+u_+)(\rho_+-\rho_-)-(\rho_+u_+-\rho_-u_-e^{-t})\nonumber\\
&=&\frac{1}{2}(\rho_++\rho_-)(u_-e^{-t}-u_+)\label{3.33},
\end{eqnarray}
from which we get (\ref{3.32}) with $w(0)=0$.

It can be derived easily from (\ref{3.32}) that
\begin{equation}\label{3.34}
\frac{d^{2}s}{dt^{2}}=-\frac{1}{2}u_{-}e^{-t}.
\end{equation}

If the delta shock exists, the following $\delta$-entropy condition
\begin{equation}\label{3.35}
u_{+}<\frac{ds}{dt}<u_{-}e^{-t},
\end{equation}
should be satisfied.

It follows from (\ref{3.22}) that
\begin{equation}\label{3.36}
\frac{ds}{dt}\mid_{t=0}=\frac{1}{2}({u_-}+u_+),
\end{equation}
which means that the speed of the delta shock waves for the inhomogeneous and the
homogeneous situations are indentical and the $\delta$-entropy condition (\ref{3.35})
is satisfied in the beginning for $u_{-}>u_{+}$ when $g(x,t,u)=-u$.

Since
\begin{equation}\label{3.37}
\frac{ds}{dt}-u_{+}=u_{-}e^{-t}-\frac{ds}{dt}=\frac{1}{2}({u_-}e^{-t}-u_+),
\end{equation}
for simplicity, we introduce the notation
\begin{equation}\label{3.38}
p(t)={u_-}e^{-t}-u_+.
\end{equation}
Differentiate (\ref{3.38}) with respect to $t$ yields
\begin{equation}\label{3.39}
p'(t)=-{u_-}e^{-t},\ \
p''(t)={u_-}e^{-t}.
\end{equation}
It follow from (\ref{3.38})
\begin{equation}\label{3.40}
p(0)={u_-}-u_+>0.
\end{equation}

In order to check the inequality (\ref{3.35}), our discussions should be
divided into four cases according to the values of $u_{-}$ and $u_{+}$ as follows.
\bigbreak

(1)
\hspace{0.05cm}
\hangafter 1
\hangindent 3.5em
\noindent
If $u_{-}>0\geq u_{+}$, then $p'(t)<0$ and $p''(t)>0$ $t\geq0$, so $p(t)$ is convex and strictly decreasing.
Moreover, $\lim\limits{t\rightarrow+\infty}p(t)=-u_+$, which implies $y=p(t)$ has the line $y=-u_+$
as its asymptotic line in the $(t,y)$-plane. So $p(t)>-u_+\geq0$ for $t\geq0$, which means that the $\delta$-entropy condition
 (\ref{3.35}) always holds for $t\geq0$ when $u_{-}>0\geq u_{+}$ and the delta shock wave always exists.

(2)\hspace{0.05cm} \hangafter 1
\hangindent 3.5em
\noindent
If $u_{-}>u_{+}>0$, it is easy to see that $p'(t)<0$ and $p''(t)>0$ for $t\geq0$. Since $p(0)>0$, there exist a unique $t_3$
 such that $p(t_3)=0$, i.e. $t_3=\ln\frac{u_-}{u_+}$. Moreover, it is easy to get that $p(t)>0$ for $0\leq t<t_3$ and $p(t)<0$ for $t>t_3$.
Thus the $\delta$-entropy condition
(\ref{3.35})holds for $0\leq t<t_3$ when $u_{-}>0\geq u_{+}$ and the delta shock wave disappear at the time $t=t_3$.

(3) \hspace{0.05cm} \hangafter 1
\hangindent 3.5em
\noindent
If $u_{-}<0$, then $p'(t)>0$ and $p''(t)<0$ $t\geq0$, so $p(t)$ is concave and strictly increasing.
So $p(t)>p(0)>0$ for $t\geq0$, which means that the $\delta$-entropy condition
(\ref{3.35}) always holds for $t\geq0$ when $u_{-}<0$ and the delta shock wave always exists.

(4)\hspace{0.05cm} \hangafter 1
\hangindent 3.5em
\noindent
If $u_{-}=0>u_+$, then $p(t)=-u_+>0$ for $t\geq0$, which means that the $\delta$-entropy condition
(\ref{3.35}) always holds for $t\geq0$ when $u_{-}=0> u_{+}$ and the delta shock wave always exists.
The proof is completed.

 \bigbreak

For the nomogeneous situation $g(x,t,u)=-u$, the disccussion about the path of the delta shock wave for
the Riemann problem (\ref{1.1}) and (\ref{1.2}) can be carried out like as before and should be divided into the following two cases.

 \bigbreak
\noindent Case 1 \hspace{0.25cm} If ${u_-}+u_+>0$, taking into account to ${u_-}>u_+$, it is easy to get $u_->0$.
So $\frac{d^{2}s}{dt^{2}}<0$ for $t\geq0$, which means that the delta shock wave curve is always concave and the speed of the
delta shock wave slows down. In the following, there are three subcases needed to be considered.

\bigbreak

(i)\hspace{0.15cm} \hangafter 1
\hangindent 3.5em
\noindent
If $u_{-}> u_{+}>0$, form the result obtainted in (2), the delta shock wave disappears at the time $t_3$ and $\frac{ds}{dt}|_{t=t_3}=u_+=u_-e^{-t_{3}}>0$.
Taking account into $\frac{d^{2}s}{dt^{2}}<0$ for $t\geq0$,
we have $\frac{ds}{dt}>\frac{ds}{dt}|_{t=t_3}>0$ for $0<t
<t_3$. So the delta shock wave curve is concave and $s(t)$
increases along with $t$ until it reaches the time $t_3$. Moreover, $\frac{ds}{dt}|_{t=t_3}=u_+=u_-e^{-t_{3}}$, so the delta shock wave curve
is tangent with the charateristic curves at the time $t_3$ on both sides of it and then disappears. After the time $t_3$ the vacuum occurs. We can
draw Fig.2(a) to depict this situation.
\bigbreak

(ii)\hspace{0.15cm} \hangafter 1
\hangindent 3.5em
\noindent
If $u_{-}> 0>u_{+}$, from the result obtained in (1), the delta shock wave always exists for $t>0$.
It is easy to get $\lim\limits{t\rightarrow+\infty}\frac{ds}{dt}=\frac{1}{2}u_+<0$. Taking account into $\frac{ds}{dt}|_{t=0}>0$ and $\frac{d^{2}s}{dt^{2}}<0$ for $t\geq0$, there exiss a unique $t_4$ such that $\frac{ds}{dt}|_{t=t_4}=0$, i.e. $t_4=\ln{-\frac{u_-}{u_+}}$. Moreover, $\frac{ds}{dt}>0$
for $0<t<t_4$ and $\frac{ds}{dt}<0$
for $t>t_4$. So the delta shock wave curve is always concave, moves forward for $0<t<t_4$, changes its direction at the time $t=t_4$ and
moves backward for $t>t_4$.  Moreover, $\lim\limits{t\rightarrow+\infty}s(t)=-\infty$, so the delta shock wave curve intersects with the $t$-axis at the time $t\hat{t}$
such that $s(\hat{t})=0$. We can
draw Fig.2(b) to depict this situation.
\bigbreak

(iii)\hspace{0.15cm} \hangafter 1
\hangindent 3.5em
\noindent
If $u_{-}> 0=u_{+}$, from the result obtained in (1), the delta shock wave always exists for $t>0$.
 Moreover, $\frac{ds}{dt}=\frac{1}{2}u_-e^{-t}>0$ and $\frac{d^{2}s}{dt^{2}}=-\frac{1}{2}u_-e^{-t}<0$ for $t\geq0$.
So the delta shock wave curve is always concave and moves forward for $t>0$. Furthermore, the delta shock wave front never intersects with the $t$-axis
since $s(t)=\frac{1}{2}u_-(1-e^{-t})>0$ for $t>0$. We can
draw Fig.2(c) to depict this situation.

 \bigbreak
\noindent Case 2 \hspace{0.25cm} If If ${u_-}+u_+\leq0$, taking into account to ${u_-}>u_+$, it is easy to get $u_+<0$.
In the following, there are three subcases needed to be considered.

\bigbreak

(i)\hspace{0.15cm} \hangafter 1
\hangindent 3.5em
\noindent
If $u_{-}> 0>u_{+}$, from the result obtained in (1), the delta shock wave always exists for $t>0$.
It is easy to get $\lim\limits{t\rightarrow+\infty}\frac{ds}{dt}=\frac{1}{2}u_+<0$. Taking account into $\frac{ds}{dt}|_{t=0}\leq0$ and $\frac{d^{2}s}{dt^{2}}<0$ for $t>0$, we have $\frac{ds}{dt}<0$
for $t>0$. So the delta shock wave curve is always concave and moves backward for $t>0$.
 Moreover, the delta shock wave curve front never intersects with the $t$-axis. We can
draw Fig.2(d) to depict this situation.
\bigbreak

(ii)\hspace{0.15cm} \hangafter 1
\hangindent 3.5em
\noindent
If $u_{+}<u_{-}<0$, from the result obtained in (3), the delta shock wave always exists for $t>0$.
It is easy to get $\lim\limits{t\rightarrow+\infty}\frac{ds}{dt}=\frac{1}{2}u_+<0$. Taking account into $\frac{ds}{dt}|_{t=0}\leq0$ and $\frac{d^{2}s}{dt^{2}}>0$ for $t\geq0$, we have $\frac{ds}{dt}<\frac{1}{2}u_+<0$ for $t>0$.
So the delta shock wave curve is always convex and moves backward for $t>0$.
 Moreover, the delta shock wave front never intersects with the $t$-axis
since $s(t)=\frac{1}{2}(u_-(1-e^{-t})+u_+t)<0$ for $t>0$. We can
draw Fig.2(e) to depict this situation.
\bigbreak

(iii)\hspace{0.15cm} \hangafter 1
\hangindent 3.5em
\noindent
If $u_{-}= 0>u_{+}$, from the result obtained in (4), the delta shock wave always exists for $t>0$. The delta shock wave curve is .
$s(t)=\frac{1}{2}u_+t$ for $t>0$. We can
draw Fig.2(f) to depict this situation.

\bigbreak
\begin{rem}\label{rem:3.2}  Similarly, we can consider the situation $g(x,t,u)=u$.
\end{rem}

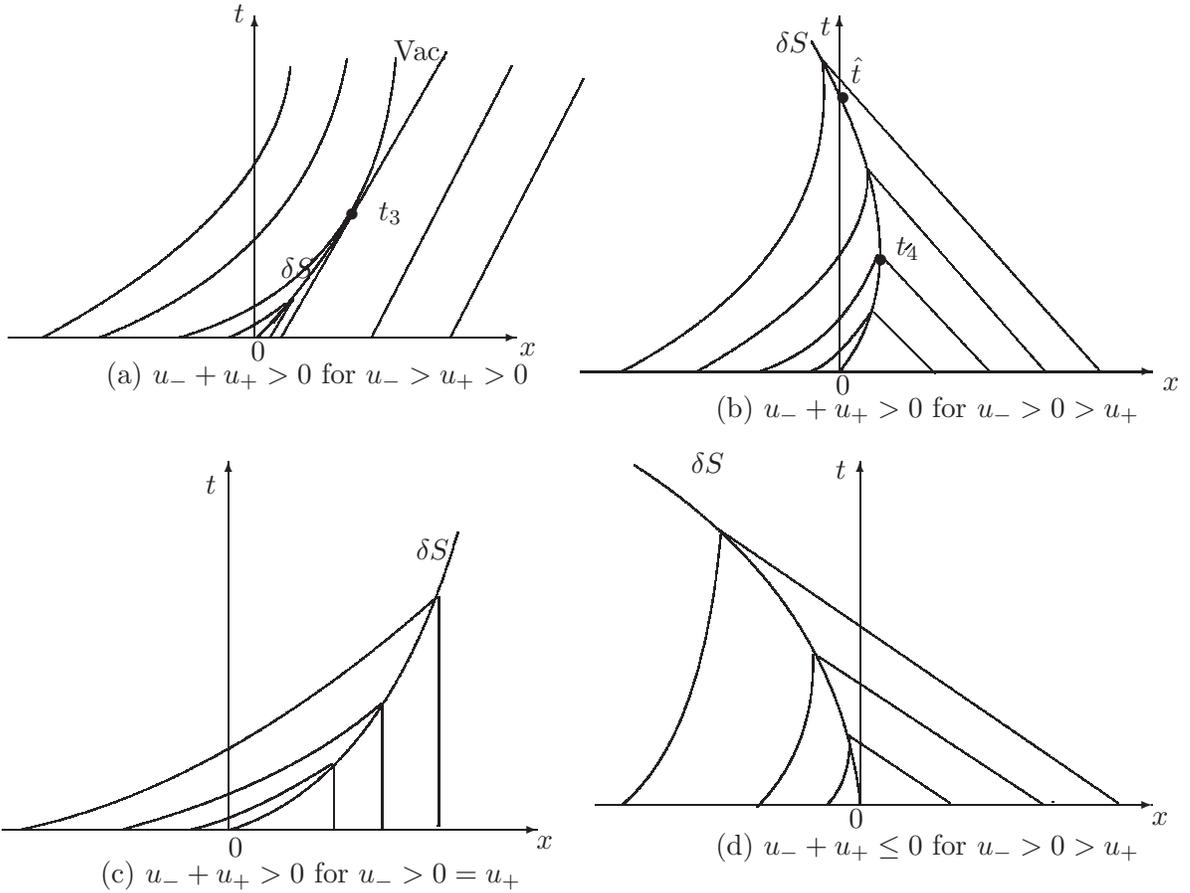
\begin{multicols}{2}
\unitlength 1mm 
\linethickness{0.4pt}
\ifx\plotpoint\undefined\newsavebox{\plotpoint}\fi 
\begin{picture}(92.25,51.25)(10,0)
\put(83.25,8.5){\vector(1,0){.07}}
\put(15.75,8.5){\line(1,0){67.5}}
\put(48.5,51.25){\vector(0,1){.07}}
\put(48.5,8.5){\line(0,1){42.75}}
\multiput(52,8.5)(.033690658,.058192956){653}{\line(0,1){.058192956}}
\multiput(64,8.5)(.033723022,.065197842){556}{\line(0,1){.065197842}}
\multiput(74.5,8.5)(.033681214,.065464896){527}{\line(0,1){.065464896}}
\qbezier(27.75,8.5)(55.75,19.5)(60.75,45.5)
\qbezier(20.25,8.5)(51.25,25.38)(53.25,44.5)
\qbezier(38.25,8.5)(64.75,16)(67.25,45.75)
\qbezier(48.75,8.5)(53.5,13)(61.25,24.5)
\qbezier(45,8.5)(48.75,10.13)(52.5,13)
\put(50.5,8.75){\line(0,1){0}}
\put(53.75,13.5){\line(-1,0){.5}}
\multiput(53.25,13.5)(-.0335821,-.0447761){67}{\line(0,-1){.0447761}}
\multiput(53.5,13.5)(-.0337079,-.0561798){89}{\line(0,-1){.0561798}}
\put(50.5,8.5){\line(0,1){0}}
\put(45.75,50.5){$t$}
\put(83.75,6.25){$x$}
\put(67,45.5){Vac.}
\put(65,24.25){$t_3$}
\put(60.5,24){$\bullet$}
\put(52,16.5){$\delta S$ }
\put(48,5.5){0}
\put(28.75,2.5){(a) $u_-+u_+>0$ for $u_->u_+>0$}
\end{picture}


\unitlength 1mm 
\linethickness{0.4pt}
\ifx\plotpoint\undefined\newsavebox{\plotpoint}\fi 
\begin{picture}(87.5,53.75)(10,0)
\put(86,6.5){\vector(1,0){.07}}
\put(10,6.5){\line(1,0){76}}
\put(44.5,53.75){\vector(0,1){.07}}
\put(44.5,6.5){\line(0,1){47.25}}
\put(68.25,48.75){\line(0,1){0}}
\put(63.75,48.25){\line(0,1){0}}
\put(73,48.5){\line(0,1){0}}
\put(70,45.75){\line(0,1){0}}
\put(77.75,45.75){\line(0,1){0}}
\put(50,6.5){\line(0,1){0}}
\put(42,51.25){$t$}
\put(87.5,4){$x$}
\put(44,3.5){0}
\put(28,0.5){(b) $u_-+u_+>0$ for $u_->0>u_+$}
\put(67.75,6.5){\line(0,1){0}}
\multiput(54,23.5)(-.032609,-.032609){23}{\line(0,-1){.032609}}
\put(82.5,6.5){\line(0,1){0}}
\qbezier(44.5,6.5)(56.88,22)(40.75,50.5)
\multiput(42,48.25)(.03372835,-.0378304467){1097}{\line(0,-1){.0378304467}}
\put(79,6.5){\line(0,1){0}}
\multiput(50,22)(.03372093,-.036046512){430}{\line(0,-1){.036046512}}
\put(64.5,6.5){\line(0,1){0}}
\put(64.5,6.5){\line(0,1){0}}
\multiput(49,14.5)(.033673469,-.033673469){245}{\line(0,-1){.033673469}}
\put(57.25,6.5){\line(0,1){0}}
\put(70.75,6.5){\line(0,1){0}}
\put(72.5,6.5){\line(0,1){0}}
\multiput(48,33.75)(.0337381916,-.0384615385){710}{\line(0,-1){.0384615385}}
\put(73,5.25){\line(0,1){0}}
\qbezier(49.5,22)(44.88,10.25)(33.75,6.5)
\qbezier(48.75,14.5)(43.88,7.5)(40.5,6.5)
\qbezier(48.25,33.25)(48.63,19.63)(25.5,6.5)
\qbezier(42.25,47.5)(44.63,21.88)(15.5,6.5)
\put(36,49){$\delta S$}
\put(52,22){$t_4$}
\put(49,20.5){$\bullet$}
\put(46,45){$\hat{t}$}
\put(44,42){$\bullet$}
\end{picture}

\end{multicols}


\begin{multicols}{2}
\unitlength 1mm 
\linethickness{0.4pt}
\ifx\plotpoint\undefined\newsavebox{\plotpoint}\fi 
\begin{picture}(87.5,55.5)(10,0)
\put(86,6.5){\vector(1,0){.07}}
\put(15,6.5){\line(1,0){71}}
\put(68.25,48.75){\line(0,1){0}}
\put(63.75,48.25){\line(0,1){0}}
\put(73,48.5){\line(0,1){0}}
\put(70,45.75){\line(0,1){0}}
\put(77.75,45.75){\line(0,1){0}}
\put(50,6.5){\line(0,1){0}}
\put(42,51.25){$t$}
\put(86,4){$x$}
\put(45,3){0}
\put(67.75,6.5){\line(0,1){0}}
\put(82.5,6.25){\line(0,1){0}}
\put(79,6.75){\line(0,1){0}}
\put(64.5,6.5){\line(0,1){0}}
\put(64.5,6.5){\line(0,1){0}}
\put(57.25,6.25){\line(0,1){0}}
\put(70.75,7){\line(0,1){0}}
\put(72.5,5.5){\line(0,1){0}}
\put(73,5.25){\line(0,1){0}}
\put(45,55.5){\vector(0,1){.07}}
\put(45,6.75){\line(0,1){48.75}}
\qbezier(45.25,6.5)(65.63,13)(75.5,46)
\put(73,37.5){\line(0,-1){30.5}}
\put(73,7){\line(-1,0){.25}}
\put(65.5,23){\line(0,-1){16.25}}
\put(65.5,6.75){\line(0,1){0}}
\put(59,15.25){\line(0,-1){8.75}}
\put(59,6.5){\line(0,1){0}}
\qbezier(72.75,37.5)(44,12)(17.25,6.5)
\qbezier(58.75,15.25)(47.25,8)(39.75,6.5)
\qbezier(65.25,23.25)(53,12.5)(30.75,6.5)
\put(70,42.5){$\delta S$}
\put(28,-0.5){(c) $u_-+u_+>0$ for $u_->0=u_+$}
\end{picture}


\unitlength 1mm 
\linethickness{0.4pt}
\ifx\plotpoint\undefined\newsavebox{\plotpoint}\fi 
\begin{picture}(87.5,52.25)(10,0)
\put(86,6.5){\vector(1,0){.07}}
\put(12,6.5){\line(1,0){73}}
\put(50,6.5){\line(0,1){0}}
\put(86,4){$x$}
\put(67.75,6.5){\line(0,1){0}}
\put(82.5,6.25){\line(0,1){0}}
\put(79,6.75){\line(0,1){0}}
\put(64.5,6.5){\line(0,1){0}}
\put(64.5,6.5){\line(0,1){0}}
\put(57.25,6.25){\line(0,1){0}}
\put(70.75,7){\line(0,1){0}}
\put(72.5,5.5){\line(0,1){0}}
\put(73,5.25){\line(0,1){0}}
\put(73,7){\line(-1,0){.25}}
\put(65.5,6.75){\line(0,1){0}}
\put(59,6.5){\line(0,1){0}}
\put(81.75,7){\line(0,1){0}}
\put(81.75,6.25){\line(0,1){0}}
\put(74,7){\line(0,1){0}}
\put(65.25,6.75){\line(0,1){0}}
\put(57.75,6.75){\line(0,1){0}}
\put(73.25,6.25){\line(0,1){0}}
\put(72.75,7){\line(0,1){0}}
\put(45.75,3.5){0}
\put(47.25,52.25){\vector(0,1){.07}}
\put(47.25,6.5){\line(0,1){45.75}}
\qbezier(47.25,6.75)(44.25,34.25)(17.25,51.75)
\multiput(28,43.5)(.0492144177,-.0337338262){1090}{\line(1,0){.0492144177}}
\put(81.25,7){\line(0,1){0}}
\multiput(41.75,26.25)(.051194539,-.033703072){586}{\line(1,0){.051194539}}
\put(71.75,6.5){\line(0,1){0}}
\multiput(45.75,15.75)(.050561798,-.033707865){267}{\line(1,0){.050561798}}
\put(59.25,6.75){\line(0,1){0}}
\qbezier(28.75,43)(26.38,15.5)(15.5,6.5)
\qbezier(41,26.5)(41.25,14.13)(33.5,6.25)
\qbezier(45.75,15.75)(46.13,10.25)(43,6.75)
\put(44,50){$t$}
\put(24.75,50.75){$\delta S$}
\put(28,0){(d) $u_-+u_+\leq0$ for $u_->0>u_+$}
\end{picture}

\end{multicols}

\newpage


\begin{multicols}{2}
\unitlength 1mm 
\linethickness{0.4pt}
\ifx\plotpoint\undefined\newsavebox{\plotpoint}\fi 
\begin{picture}(87.5,55.5)(10,0)
\put(86,6.5){\vector(1,0){.07}}
\put(15,6.5){\line(1,0){71}}
\put(50,6.5){\line(0,1){0}}
\put(85.5,4){$x$}
\put(67.75,6.5){\line(0,1){0}}
\put(82.5,6.25){\line(0,1){0}}
\put(79,6.75){\line(0,1){0}}
\put(64.5,6.5){\line(0,1){0}}
\put(64.5,6.5){\line(0,1){0}}
\put(57.25,6.25){\line(0,1){0}}
\put(70.75,7){\line(0,1){0}}
\put(72.5,5.5){\line(0,1){0}}
\put(73,5.25){\line(0,1){0}}
\put(73,7){\line(-1,0){.25}}
\put(65.5,6.75){\line(0,1){0}}
\put(59,6.5){\line(0,1){0}}
\put(81.75,7){\line(0,1){0}}
\put(81.75,6.25){\line(0,1){0}}
\put(74,7){\line(0,1){0}}
\put(65.25,6.75){\line(0,1){0}}
\put(57.75,6.75){\line(0,1){0}}
\put(73.25,6.25){\line(0,1){0}}
\put(72.75,7){\line(0,1){0}}
\put(50,3.5){0}
\put(81.25,7){\line(0,1){0}}
\put(71.75,6.5){\line(0,1){0}}
\put(59.25,6.75){\line(0,1){0}}
\put(48,52){$t$}
\put(16,40){$\delta S$}
\put(51,55.5){\vector(0,1){.07}}
\put(51,6.5){\line(0,1){48.75}}
\put(51,55.5){\vector(0,1){.07}}
\put(51,55.5){\line(0,1){0}}
\qbezier(51,6.5)(28,11.75)(13,41.5)
\qbezier(15.5,36)(17.38,17.13)(23.75,6.5)
\multiput(16,36.25)(.074106113,-.0337370242){880}{\line(1,0){.074106113}}
\multiput(25.5,22.5)(.079473684,-.033684211){475}{\line(1,0){.079473684}}
\put(63.25,6.5){\line(0,1){0}}
\put(54.75,6.5){\line(0,1){0}}
\multiput(34.75,14.25)(.087755102,-.033673469){230}{\line(1,0){.087755102}}
\put(56.25,6){\line(0,1){0}}
\qbezier(25.5,22.5)(25.63,14.38)(30.25,6.5)
\qbezier(34.25,14.5)(35.13,8.38)(37.5,6.5)
\put(31,1){(e) $u_-+u_+\leq0$ for $u_+<u_-<0$}
\end{picture}


\unitlength 1mm 
\linethickness{0.4pt}
\ifx\plotpoint\undefined\newsavebox{\plotpoint}\fi 
\begin{picture}(87.5,55.5)(10,-2)
\put(86,6.5){\vector(1,0){.07}}
\put(15,6.5){\line(1,0){71}}
\put(50,6.5){\line(0,1){0}}
\put(86,4){$x$}
\put(67.75,6.5){\line(0,1){0}}
\put(82.5,6.25){\line(0,1){0}}
\put(79,6.75){\line(0,1){0}}
\put(64.5,6.5){\line(0,1){0}}
\put(64.5,6.5){\line(0,1){0}}
\put(57.25,6.25){\line(0,1){0}}
\put(70.75,7){\line(0,1){0}}
\put(72.5,5.5){\line(0,1){0}}
\put(73,5.25){\line(0,1){0}}
\put(73,7){\line(-1,0){.25}}
\put(65.5,6.75){\line(0,1){0}}
\put(59,6.5){\line(0,1){0}}
\put(81.75,7){\line(0,1){0}}
\put(81.75,6.25){\line(0,1){0}}
\put(74,7){\line(0,1){0}}
\put(65.25,6.75){\line(0,1){0}}
\put(57.75,6.75){\line(0,1){0}}
\put(73.25,6.25){\line(0,1){0}}
\put(72.75,7){\line(0,1){0}}
\put(49.5,3.5){0}
\put(81.25,7){\line(0,1){0}}
\put(71.75,6.5){\line(0,1){0}}
\put(59.25,6.75){\line(0,1){0}}
\put(48,52.5){$t$}
\put(25.75,49){$\delta S$}
\put(51,55.5){\vector(0,1){.07}}
\put(51,6.75){\line(0,1){48.75}}
\put(51,55.5){\vector(0,1){.07}}
\put(51,55.5){\line(0,1){0}}
\put(63.25,6.5){\line(0,1){0}}
\put(54.75,6.5){\line(0,1){0}}
\put(56.25,6){\line(0,1){0}}
\multiput(51,6.5)(-.033707865,.064606742){623}{\line(0,1){.064606742}}
\multiput(33,41)(.0444770283,-.0337243402){1023}{\line(1,0){.0444770283}}
\multiput(39.5,28.75)(.0425925926,-.0337037037){660}{\line(1,0){.0425925926}}
\put(68.25,6){\line(0,1){0}}
\multiput(45.75,16.75)(.042813456,-.033639144){310}{\line(1,0){.042813456}}
\multiput(33.25,40)(-.03373494,-.080722892){415}{\line(0,-1){.080722892}}
\put(19.25,6.5){\line(0,1){0}}
\multiput(39.5,28.75)(-.033673469,-.091836735){245}{\line(0,-1){.091836735}}
\put(31.25,6.25){\line(0,1){0}}
\multiput(45.75,16.5)(-.03365385,-.09615385){104}{\line(0,-1){.09615385}}
\put(42.25,6.5){\line(0,1){0}}
\put(31,0.5){(f) $u_-+u_+\leq0$ for $u_-=0>u_+$}
\end{picture}
\put(-78,-3){\makebox(0,0)[cc]
{Fig.2 The delta shock wave solution to (1.1) and (1.2) when $f(x,t,u)=0$ and $g(x,t,u)=-u$. }}

\end{multicols}

\section{The situation for $f(x,t,u)\neq0$ and $g(x,t,u)\neq0$}
In this section. we consider the case that $f(x,t,u)\neq0$, $g(x,t,u)\neq0$ and $f(x,t,u)\neq g(x,t,u)$. For simplicity, we take $f(x,t,u)=1$ and
$g(x,t,u)=-u$ to display how the delta shock front develop under the effect of the source term.
Similar to subsection 3.1,
we have the following theorem to depict the delta shock wave solution to the Riemann problem (\ref{1.1}) and
(\ref{1.2}) when $u_{-}>u_{+}$.

\begin{thm}\label{thm:4.1}
For the situation $f(x,t,u)=1$ and $g(x,t,u)=-u$, if $u_{-}>u_{+}$,
then the Riemann solution to (\ref{1.1}) and (\ref{1.2}) is a delta
shock solution which can be express as
\begin{equation}\label{4.1}
(\rho,u)(x,t)=\left\{\begin{array}{ll}
(\rho_-,u_-e^{-t}),\ \ x<s(t),\\
(w(t),u_{\delta}(t)),\ \ x=s(t)\\
(\rho_+,u_+t),\ \ x>s(t),
\end{array}\right.
\end{equation}
with
\begin{equation}\label{4.2}
 s(t)=\frac{1}{2}(u_{-}(1-e^{-t})+u_{+}t+\frac{1}{2}t^{2}),
\end{equation}
\begin{equation}\label{4.3}
 w(t)=\frac{1}{2}(\rho_{-}+\rho_{+})(u_{-}(1-e^{-t})-u_{+}t-\frac{1}{2}t^{2}),
\end{equation}
\begin{equation}\label{4.4}
 u_{\delta}(t)=\frac{1}{2}(u_{-}e^{-t}+u_{+}+t).
\end{equation}
Furthermore, the delta shock disappears in finite time.
\end{thm}

For $x<s(t)$, similar to Section 3.2, we have
\begin{equation}\label{4.5}
x=u_-(1-e^{-t})+x_0, for t\geq0,x_{0}<0.
\end{equation}
and
\begin{equation}\label{4.6}
(\rho,u)(x,t)=(\rho_-,u_-e^{-t}), for all \ \ x<s(t), t\geq0.
\end{equation}

For $x>s(t)$, similar to Section 3.1, we have
\begin{equation}\label{4.7}
x=\frac{1}{2}t^{2}+u_+t+x_0, for t\geq0,x_{0}>0.
\end{equation}
and
\begin{equation}\label{4.8}
(\rho,u)(x,t)=(\rho_+,u_++t), for all \ \ x>s(t), t\geq0.
\end{equation}

It is clear to see that after some time $t$ the delta shock wave will disappear for the reason that
the $\delta$-entropy condition cannot be satisfied. However, for a sufficiently small time $t$, the delta shock wave may exist and should satisfy
the generalized Rankine-Hongniot condition (\ref{1.5}), so we have
\begin{equation}\label{4.9}
\frac{ds}{dt}(u_++t-u_-e^{-t} )=\frac{1}{2}({u_++t}^2-(u_-e^{-t})^2),
\end{equation}
namely, we have
\begin{equation}\label{4.10}
u_{\delta}(t)=\sigma(t)=\frac{ds}{dt}=\frac{1}{2}({u_-}e^{-t}+u_++t).
\end{equation}
So (\ref{4.4})is obtained. With $s(0)=0$ in mind, we can get the expression of the delta shock
front (\ref{4.2}) in finite time.

Moreover, the weight of the delta shock can be get from the second
equality of (\ref{1.5}) with
\begin{eqnarray}
\frac{d
w(t)}{dt}&=&\sigma(t)(\rho_+-\rho_-)-(\rho_+(u_++t)-\rho_-u_-e^{-t})\nonumber\\
&=&\frac{1}{2}({u_-}e^{-t}+u_++t)(\rho_+-\rho_-)-(\rho_+(u_++t)-\rho_-u_-e^{-t})\nonumber\\
&=&\frac{1}{2}(\rho_++\rho_-)(u_-e^{-t}-(u_++t))\label{4.11},
\end{eqnarray}
from which we get (\ref{4.3}) with $w(0)=0$.

It can be derived easily from (\ref{4.10}) that
\begin{equation}\label{4.12}
\frac{d^{2}s}{dt^{2}}=\frac{1}{2}(1-u_{-}e^{-t}).
\end{equation}

If the delta shock exists, the following $\delta$-entropy condition
\begin{equation}\label{4.13}
u_{+}+t<\frac{ds}{dt}<u_{-}e^{-t},
\end{equation}
should be satisfied.

It follows from (\ref{4.10}) that
\begin{equation}\label{4.14}
\frac{ds}{dt}\mid_{t=0}=\frac{1}{2}({u_-}+u_+),
\end{equation}
which means that the speed of the delta shock waves for the inhomogeneous and the
homogeneous situations are indentical and the $\delta$-entropy condition (\ref{4.13})
is satisfied in the beginning for $u_{-}>u_{+}$ when $f(x,t,u)=1$ and $g(x,t,u)=-u$.

Since
\begin{equation}\label{4.15}
\frac{ds}{dt}-(u_{+}+t)=u_{-}e^{-t}-\frac{ds}{dt}=\frac{1}{2}({u_-}e^{-t}-(u_++t)),
\end{equation}
for simplicity, we introduce the notation
\begin{equation}\label{4.16}
q(t)={u_-}e^{-t}-(u_++t).
\end{equation}

Differentiate (\ref{4.16}) with respect to $t$ yields
\begin{equation}\label{4.17}
q'(t)=-{u_-}e^{-t}-1,\ \
q''(t)={u_-}e^{-t}.
\end{equation}
It follow from (\ref{4.16})
\begin{equation}\label{4.18}
q(0)={u_-}-u_+>0.
\end{equation}
Since $\lim\limits{t\rightarrow+\infty}q(t)=-\infty$, there exists $t_5>0$ such that $q(t_5)=0$, i.e., $u_-e^{-t_5}-(u_++t_5)=0$.
Next, we will prove that $t_5$ is unique and $q(t)>0$ for $0< t<t_5$.
our discussions should be
divided into two cases according to the values of $u_{-}$ and $u_{+}$ as follows.
\bigbreak

(1)
\hspace{0.05cm}
\hangafter 1
\hangindent 3.5em
\noindent
If $u_{-}\geq0$, then $q'(t)<0$, which means that $q(t)$ is strictly decreasing for $t\geq0$. Obviouly, $t_5$ is unique and $q(t)>q(t_5)=0$ for $0< t<t_5$.

(2)\hspace{0.05cm} \hangafter 1
\hangindent 3.5em
\noindent
If $u_-<0$, since $q'(0)=-u_--1$, then it should be divided into the following two subcases.

\hspace{0.3cm} (2a)\hspace{0.15cm} \hangafter 1
\hangindent 5em
\noindent
\noindent
If $-1\leq u_{-}<0$, then $q'(0)\leq0$. Since $p''(t)<0$ for $t\geq0$, so $q'(t)<q'(0)\leq0$, which means that $q(t)$ is strictly decreasing for $t\geq0$.
Obviouly, $t_5$ is unique and $q(t)>q(t_5)=0$ for $0< t<t_5$.

\hspace{0.3cm} (2b)\hspace{0.15cm} \hangafter 1
\hangindent 5em
\noindent
\noindent
If $ u_{-}<-1$, then $q'(0)>0$. Since $q''(t)<0$ for $t\geq0$ and $\lim\limits{t\rightarrow+\infty}q'(t)=-1<0$, there exists $\tilde{t}>0$ such that $q'(\tilde{t})=0$.
Moreover, $q'(t)>q'(0)>0$ for $0< t<\tilde{t}$ and $q'(t)<q'(\tilde{t})=0$ for $t>\tilde{t}$. Since $q(0)>0$, we have $q(t)>0$ for $0< t\leq\tilde{t}$.
Since $q(t_5)=0$, we conclude that $t_5>\tilde{t}$.
Combing with $q'(t)<0$ for $t>\tilde{t}$, which means that $q(t)$ is strictly decreasing for $t>\tilde{t}$, we can claim that $t_5$ is unique and $q(t)>q(t_5)=0$ for $\tilde{t}< t<t_5$. In a word, $q(t)>0$ for $0< t<t_5$.

\bigbreak
From the above discussion, we can conclude that there exists a unique $t_5>0$ such that $q(t_5)=0$ and $q(t)>0$ for $0< t<t_5$, which means that
the the $\delta$-entropy condition
(\ref{4.13}) holds for $0\leq t<t_5$. Moreover, $\frac{ds}{dt}|_{t=t_5}=u_-e^{-t_5}=u_++t_5$,
So the delta shock wave is tangent with the characteristic curve at the time $t_5$ on both sides of it and then disappear.
The proof is completed.

\bigbreak

For the nomogeneous situation $f(x,t,u)=1$ and $g(x,t,u)=-u$, from Theorem 4.1, the delta shock wave disappears at the time $t_5$. Next the disccussion about the path of the delta shock wave for the Riemann problem (\ref{1.1}) and (\ref{1.2}) can be carried out like as before.
Since $\frac{ds}{dt}\mid_{t=0}=\frac{1}{2}({u_-}+u_+)$, then it should be divided into the following three cases. Hereafter, for $u_+<0$, we denote $t_2^\ast=-u_+$ the symmetry axises of the characteristics on the right hand side of the delta shock wave.

 \bigbreak
\noindent Case 1 \hspace{0.25cm} If ${u_-}+u_+> 0$, taking into account to ${u_-}>u_+$, it is easy to get $u_->0$. So $\frac{ds}{dt}|_{t=t_5}=u_-e^{-t_5}=u_++t_5>0$, from which we have $t_5-t_2^\ast=u_-e^{-t_5}>0$, i.e., $t_5>t_2^\ast$. From (\ref{4.12}),
we have $\frac{d^{3}s}{dt^{3}}=\frac{1}{2}u_{-}e^{-t}>0$ for $t\geq0$ and $\frac{d^{2}s}{dt^{2}}|_{t=0}=\frac{1}{2}(1-u_-)$.
In the following, there are two subcases needed to be considered.
\bigbreak

(i)\hspace{0.15cm} \hangafter 1
\hangindent 3.5em
\noindent
If $0<u_{-}\leq1$, then $\frac{d^{2}s}{dt^{2}}|_{t=0}=\frac{1}{2}(1-u_-)>0$. Since $\frac{d^{3}s}{dt^{3}}=\frac{1}{2}u_{-}e^{-t}>0$,
$\frac{d^{2}s}{dt^{2}}>0$ for $t>0$. Taking into account to $\frac{ds}{dt}\mid_{t=0}=\frac{1}{2}({u_-}+u_+)>0$, we have $\frac{ds}{dt}>0$ for $t>0$.
Since the delta shock wave disappear at the time $t_5$, the delta shock wave is convex and increases along with $t$ such that the delta shock wave move forward for $0\leq t<t_5$ and then disappears. We draw Fig.3(a) and Fig.3(b) to depict this situation according to $u_+\geq0$ and $u_+<0$, where $t=t_2^\ast=-u_+$ is the the symmetry axis of the characteristic
curves on the right hand side of the delta shock waves and $t_2^\ast<t_5$.

\bigbreak

(ii)\hspace{0.15cm} \hangafter 1
\hangindent 3.5em
\noindent
If $u_{-}> 1$, then $\frac{d^{2}s}{dt^{2}}|_{t=0}=\frac{1}{2}(1-u_-)<0$. Since $\frac{d^{3}s}{dt^{3}}=\frac{1}{2}u_{-}e^{-t}>0$, and
$\lim\limits{t\rightarrow+\infty}\frac{d^{2}s}{dt^{2}}=\frac{1}{2}>0$, there exists a unique $t_6$ such that $\frac{d^{2}s}{dt^{2}}|_{t=t_6}=0$,
i.e. $t_6=\ln u_-$. Moreover, $\frac{d^{2}s}{dt^{2}}<0$ for $0\leq t<t_6$ and $\frac{d^{2}s}{dt^{2}}>0$ for $t>t_6$.
At the time $t_6$, $\frac{ds}{dt}|_{t=t_6}=\frac{1}{2}(u_-e^{-t_6}+u_++t_6)=\frac{1}{2}(1+u_++\ln u_-)$. Then, it should be divided into the following three subcases.

\bigbreak
\hspace{0.4cm} (a)\hspace{0.15cm} \hangafter 1
\hangindent 5em
\noindent
\noindent
If $1+u_++\ln u_-=0$, which is equivalent to $u_-e^{-t_6}+u_++t_6=0$. Furthermore, it is easy to get $u_+<0$. Since $q(t_5)=u_-e^{-t_5}-(u_++t_5)=0$, we have $t_5-t_6=u_-(e^{-t_5}+e^{-t_6})>0$,
i.e. $t_5>t_6$.
which mean that the delta shock wave disapppears after the time $t_6$. Moreover, $\frac{ds}{dt}>0$ for $0< t<t_6$ and $t_6<t<t_5$, $\frac{d^{2}s}{dt^{2}}<0$ for $0\leq t<t_6$ and $\frac{d^{2}s}{dt^{2}}>0$ for $t_6<t<t_5$, which means that the delta shock wave is concave for $0\leq t<t_6$, convex for $t_6<t<t_5$, always increases along with $t$ such that the delta shock wave shoud always move forward for $0\leq t<t_5$ and then disappears. We draw Fig.3(c) to depict this situation, where $t_6<t_2^\ast<t_5$.

\hspace{0.4cm} (b)\hspace{0.15cm} \hangafter 1
\hangindent 5em
\noindent
If $1+u_++\ln u_-<0$, which is equivalent to $u_-e^{-t_6}+u_++t_6<0$. Furthermore, it is easy to get $u_+<0$. Since $q(t_5)=u_-e^{-t_5}-(u_++t_5)=0$, we have $t_5-t_6>u_-(e^{-t_5}+e^{-t_6})>0$,
i.e. $t_5>t_6$.
which mean that the delta shock wave disapppears after the time $t_6$. Moreover, $\frac{ds}{dt}|_{t=0}>0$, $\frac{ds}{dt}|_{t=t_6}<0$, $\frac{ds}{dt}|_{t=t_5}=>0$,
combining with $\frac{d^{2}s}{dt^{2}}<0$ for $0\leq t<t_6$ and $\frac{d^{2}s}{dt^{2}}>0$ for $t_6<t<t_5$, there exists a unique $t_7$ with
$t_7<t_6$ and a unique $t_8$ with $t_6<t_8<t_5$
such that $\frac{ds}{dt}|_{t=t_7}=\frac{ds}{dt}|_{t=t_8}=0$. So $\frac{ds}{dt}>0$ for $0< t<t_7$ and $t_8<t<t_5$ and $\frac{ds}{dt}<0$ for $t_7<t<t_8$.
In a word, the delta shock wave is concave for $0\leq t<t_6$, convex for $t_6<t<t_5$. Furthermore, the delta shock wave increases along with $t$ such that the delta shock wave move forward for $0\leq t<t_7$ and $t_8<t<t_5$, decreases along with $t$ such that the delta shock wave move backward for $t_7<t<t_8$ and then disappears.
We draw Fig.3(d) to depict this situation, where $t_8<t_2^\ast<t_5$.

\hspace{0.4cm} (c)\hspace{0.15cm} \hangafter 1
\hangindent 5em
\noindent
\noindent
If $1+u_++\ln u_->0$, which is equivalent to $u_-e^{-t_6}+u_++t_6>0$ and $u_+>-(1+\ln u_-)$. If $t_5>t_6$, the delta shock wave is similar to Case 1 (ii)(a) in
this section, see Fig.3(c). If $t_5\leq t_6$, $\frac{ds}{dt}>0$ and $\frac{d^{2}s}{dt^{2}}<0$ for $0< t<t_5$, In a word, the delta shock wave is concave, increases along with $t$ such that the delta shock wave move forward for $0\leq t<t_5$, and then disappears.
We draw Fig.3(e) and Fig.3(f) to depict this situation according to $u_+\geq0$ and $-(1+\ln u_-)<u_+<0$, where $t_2^\ast<t_5$.

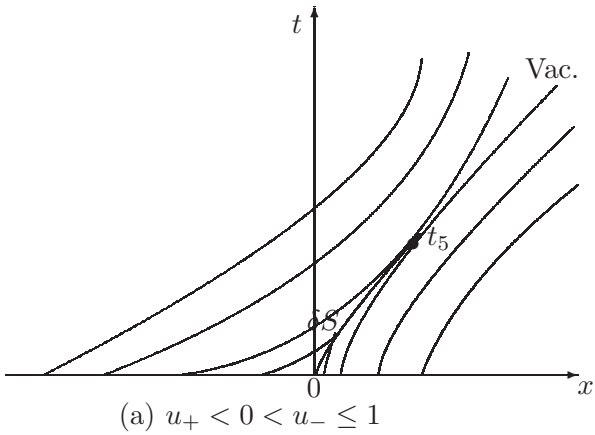
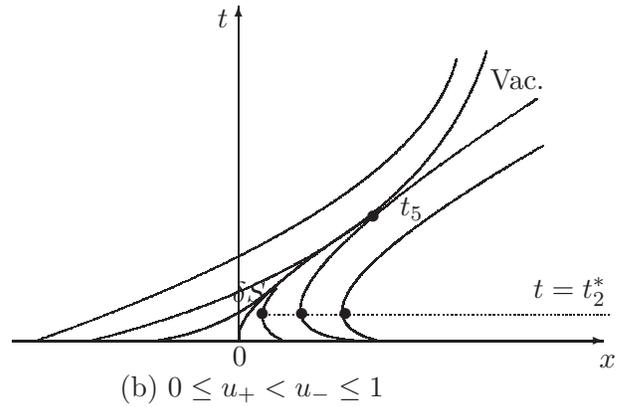
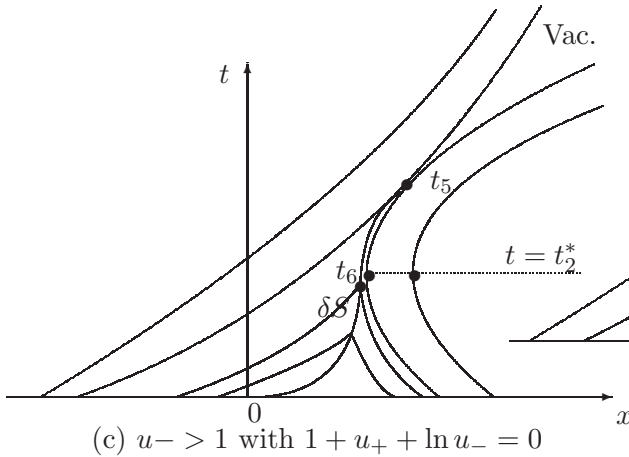
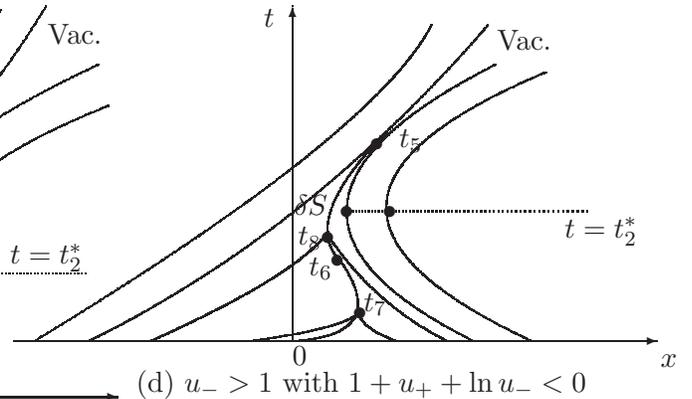
\begin{multicols}{2}
\unitlength 1mm 
\linethickness{0.4pt}
\ifx\plotpoint\undefined\newsavebox{\plotpoint}\fi 
\begin{picture}(87.5,55.5)(10,0)
\put(86,6.5){\vector(1,0){.07}}
\put(10,6.5){\line(1,0){76}}
\put(50,6.5){\line(0,1){0}}
\put(86,4){$x$}
\put(73,5.25){\line(0,1){0}}
\put(50,3.5){0}
\put(48,52){$t$}
\put(51,55.5){\vector(0,1){.07}}
\put(51,6.5){\line(0,1){48.75}}
\put(51,55.5){\vector(0,1){.07}}
\put(51,55.5){\line(0,1){0}}
\put(19.25,6.5){\line(0,1){0}}
\put(31.25,6.5){\line(0,1){0}}
\put(42.25,6.5){\line(0,1){0}}
\qbezier(23,6.5)(64.25,22.5)(71.5,49.25)
\qbezier(54.5,6.5)(55.63,16.75)(83.25,45)
\qbezier(53.75,12)(52.5,8.13)(52.25,6.5)
\qbezier(59.5,6.5)(60.5,15.5)(85.5,39.5)
\qbezier(65,25.25)(51.88,10.5)(51.25,6.5)
\qbezier(33.25,6.5)(61,10.25)(76.75,46)
\qbezier(44,6.5)(49.88,8.25)(54.25,12)
\qbezier(15,6.5)(63.88,31.5)(65.25,48.5)
\qbezier(65.25,6.5)(69.88,18.75)(86,31.75)
\put(66,23.75){$t_5$}
\put(63.15,23){$\bullet$}
\put(79,46){Vac.}
\put(50,12.5){$\delta S$}
\put(25,0){(a) $u_+<0<u_-\leq1$}
\end{picture}

\unitlength 1mm 
\linethickness{0.4pt}
\ifx\plotpoint\undefined\newsavebox{\plotpoint}\fi 
\begin{picture}(90.25,51.75)(10,0)
\put(89.25,7.25){\vector(1,0){.07}}
\put(10.75,7.25){\line(1,0){78.5}}
\put(40.75,51.75){\vector(0,1){.07}}
\put(40.75,7.25){\line(0,1){44.25}}
\qbezier(21,7.25)(65.13,18.38)(73.75,45.75)
\qbezier(57.75,7.25)(32.38,8.13)(80.5,39.5)
\qbezier(41,7.25)(39.13,11.5)(59.75,24.75)
\qbezier(45.75,14.25)(41.75,9.63)(46.75,7.25)
\qbezier(29.5,7.25)(39,8.63)(45.5,14.5)
\qbezier(59.5,7.25)(43.13,11.13)(81.25,33.25)
\bezier{80}(44.25,10.75)(70,10.75)(90,10.75)
\qbezier(13.75,7.25)(64.75,25.38)(69.75,44.75)
\put(38,49){$t$}
\put(88.75,3.75){$x$}
\put(62.25,24){$t_5$}
\put(57.75,22.95){$\bullet$}
\put(74.25,40.75){Vac.}
\put(80,13){$t=t_2^\ast$}
\put(48.25,10){$\bullet$}
\put(43,10){$\bullet$}
\put(54,10){$\bullet$}
\put(40,4){0}
\put(40,12.5){$\delta S$}
\put(25,0){(b) $0\leq u_+<u_-\leq1$}
\end{picture}

\end{multicols}

\begin{multicols}{2}
\unitlength 1mm 
\linethickness{0.4pt}
\ifx\plotpoint\undefined\newsavebox{\plotpoint}\fi 
\begin{picture}(89.75,59.25)(0,0)
\put(89.25,7.25){\vector(1,0){.07}}
\put(8.75,7.25){\line(1,0){80.5}}
\put(40.75,51.75){\vector(0,1){.07}}
\put(40.75,7.25){\line(0,1){44.25}}
\put(37,49){$t$}
\put(89.75,3.75){$x$}
\put(40.75,4){0}
\qbezier(55.75,22.25)(55,7.25)(41.25,7.25)
\qbezier(61,34.5)(56.23,31.43)(55.75,22.5)
\qbezier(18,7.25)(56.13,20.75)(79.75,59.25)
\qbezier(13.25,7.25)(56.5,33.38)(73.75,58.75)
\qbezier(66.25,7.25)(39.5,29.63)(86.75,51.5)
\qbezier(73.5,7.25)(46.25,29.38)(88,46)
\qbezier(31.25,7.25)(48.5,13.13)(55.75,22.25)
\qbezier(64,7.25)(56.63,12.75)(55.75,22)
\qbezier(36.75,7.25)(50.88,12.38)(54.5,15.75)
\qbezier(54.75,15.25)(57.88,7.88)(60.5,7.25)
\put(52.5,23){$t_6$}
\put(56,22.5){$\bullet$}
\put(62,22.5){$\bullet$}
\put(54.85,21){$\bullet$}
\put(65,35){$t_5$}
\put(61,34.5){$\bullet$}
\bezier{60}(56.5,23.75)(70.5,23.75)(85,23.75)
\put(75,25){$t=t_2^\ast$}
\put(50,18){$\delta S$}
\put(20,0.5){(c)  $u->1$
with $1+u_++\ln {u_-}=0$}
\put(80,54.25){Vac.}
\end{picture}

\unitlength 1mm 
\linethickness{0.4pt}
\ifx\plotpoint\undefined\newsavebox{\plotpoint}\fi 
\begin{picture}(89.75,51.75)(10,0)
\put(89.25,7.25){\vector(1,0){.07}}
\put(3.75,7.25){\line(1,0){85.5}}
\put(40.75,51.75){\vector(0,1){.07}}
\put(40.75,7.25){\line(0,1){44.25}}
\put(37,49){$t$}
\put(89.75,3.75){$x$}
\put(40.75,4){0}
\qbezier(66.75,49.25)(42.75,21.5)(13.75,7.25)
\qbezier(59.25,49.25)(51.63,34.63)(6.5,7.25)
\qbezier(64.75,7.25)(29.75,24.75)(67.75,44)
\qbezier(74.5,43)(33,25.63)(72.5,7.25)
\bezier{60}(48,24.5)(60,24.5)(80,24.5)
\qbezier(47.25,17.75)(53.75,8.13)(41.25,7.25)
\qbezier(53.75,35.5)(41.5,23.38)(47.25,17.75)
\qbezier(45.5,21.75)(41.25,15.88)(22,7.25)
\qbezier(45.5,21.25)(51.63,12.13)(61.25,7.25)
\qbezier(49.5,11)(50,8.88)(54.5,7.25)
\qbezier(49.5,11)(48,9)(35.5,7.25)
\put(68,46){Vac.}
\put(55,33){$t_5$}
\put(51,32.5){$\bullet$}
\put(41.5,20){$t_8$}
\put(44.5,20){$\bullet$}
\put(43,16){$t_6$}
\put(52.75,23.5){$\bullet$}
\put(47,23.5){$\bullet$}
\put(45.75,17){$\bullet$}
\put(50.25,11.25){$t_7$}
\put(48.75,10){$\bullet$}
\put(77,21){$t=t_2^\ast$}
\put(41,24){$\delta S$}
\put(20,0.5){(d) $u_->1$
with $1+u_++\ln {u_-}<0$}
\end{picture}

\end{multicols}

\begin{multicols}{2}
\unitlength 1mm 
\linethickness{0.4pt}
\ifx\plotpoint\undefined\newsavebox{\plotpoint}\fi 
\begin{picture}(89.75,51.75)(10,0)
\put(89.25,7.25){\vector(1,0){.07}}
\put(13.75,7.25){\line(1,0){75.5}}
\put(40.75,51.75){\vector(0,1){.07}}
\put(40.75,7.25){\line(0,1){44.25}}
\put(37,49){$t$}
\put(89.75,3.75){$x$}
\put(40.75,4){0}
\qbezier(54.25,7.25)(57,26.5)(79.75,35.75)
\qbezier(61.75,7.25)(61,25.25)(85.25,33.25)
\qbezier(40.75,7.25)(46.25,6.75)(62.75,25)
\qbezier(52,13.5)(50.13,11.75)(49.75,7.25)
\qbezier(51.75,14.25)(42.88,9.13)(36.5,7.25)
\qbezier(18.5,7.25)(57.75,17.13)(69,42.5)
\qbezier(29.75,7.25)(65.75,19.63)(73.75,40.25)
\put(77,37.75){Vac.}
\put(66.25,25.75){$t_5$}
\put(63,25){$\bullet$}
\put(49.5,15){$\delta S$}
\put(18.25,1){(e)  $u_->1$
with $1+u_++\ln {u_-}>0$, }
\put(30,-2.5){ $t_5\leq t_6$ and $u_+\geq0$}
\end{picture}


\unitlength 1mm 
\linethickness{0.4pt}
\ifx\plotpoint\undefined\newsavebox{\plotpoint}\fi 
\begin{picture}(89.75,51.75)(10,0)
\put(89.25,7.25){\vector(1,0){.07}}
\put(13.75,7.25){\line(1,0){75.5}}
\put(40.75,51.75){\vector(0,1){.07}}
\put(40.75,7.25){\line(0,1){44.25}}
\put(37,49){$t$}
\put(89.75,3.75){$x$}
\put(40.75,4){0}
\qbezier(40.75,7.25)(46.25,6.75)(62.75,25)
\qbezier(51.75,14.25)(42.88,9.13)(36.5,7.25)
\qbezier(18.5,7.25)(57.75,17.13)(69,42.5)
\qbezier(29.75,7.25)(65.75,19.63)(73.75,40.25)
\put(77,37.75){Vac.}
\put(66.25,25.75){$t_5$}
\put(63,25){$\bullet$}
\put(60,19.5){$\bullet$}
\put(65,19.5){$\bullet$}
\qbezier(80.5,7.25)(40.38,20.88)(82.75,34.75)
\qbezier(86.25,32.25)(45.5,22.63)(85.75,7.25)
\bezier{50}(61.25,20.75)(70,20.75)(88,20.75)
\put(87.75,20.75){\line(0,1){0}}
\qbezier(54.25,16)(54.63,10.75)(64.5,7.5)
\put(84.75,17){$t=t_2^\ast$}
\put(20,12.5){$\delta S$}
\put(18.25,1){(f) $u_->1$
with $1+u_++\ln {u_-}>0$,}
\put(30,-2.5){$t_5\leq t_6$ and $u_+<0$}
\end{picture}
\put(-80,-6){\makebox(0,0)[cc]
{Fig.3 The delta shock wave solution to (1.1) and (1.2) when $f(x,t,u)=1$ }}
\put(-80,-10){\makebox(0,0)[cc]
{and $g(x,t,u)=-u$ for $u_-+u_+>0$}}
\end{multicols}

 \bigbreak
\noindent Case 2 \hspace{0.25cm} If ${u_-}+u_+<0$, taking into account to ${u_-}>u_+$, it is easy to get $u_+<0$. According to (\ref{4.12}),
$\frac{d^{2}s}{dt^{2}}|_{t=0}=\frac{1}{2}(1-u_-)$, there are two subcases needed to be considered as follows.
\bigbreak

(i)\hspace{0.15cm} \hangafter 1
\hangindent 3.5em
\noindent
If $u_{-}\leq1$, then $\frac{d^{2}s}{dt^{2}}|_{t=0}=\frac{1}{2}(1-u_-)\geq0$.
we claim that $\frac{d^{2}s}{dt^{2}}>0$ for $t>0$. In fact, the discussions should be divided into the following two subcases
according to $u_-\leq0$ and $0<u_-\leq1$.

\bigbreak
\hspace{0.4cm} (a)\hspace{0.15cm} \hangafter 1
\hangindent 5em
\noindent
\noindent
If $u_{-}\leq0$,
it is obvious that $\frac{d^{2}s}{dt^{2}}=\frac{1}{2}(1-u_-e^{-t})>0$ for $t>0$. Since at the time $t_5$ when the delta shock wave disappears $\frac{ds}{dt}|_{t=t_5}=u_-e^{-t_5}\leq0$, combining with $\frac{ds}{dt}|_{t=0}=\frac{1}{2}(u_-+u_+)<0$, we have $\frac{ds}{dt}<0$ for $0<t<t_5$. So the delta shock wave is convex for $0\leq t<t_5$, decreases along with $t$ such that the delta shock wave move fbackward for $0\leq t<t_5$ and then disappear. We draw Fig.4(a) to depict this situation, where $t_5\leq t_2^\ast$.

\hspace{0.4cm} (b)\hspace{0.15cm} \hangafter 1
\hangindent 5em
\noindent
\noindent
If $0<u_{-}\leq1$, since $\frac{d^{3}s}{dt^{3}}=\frac{1}{2}u_{-}e^{-t}>0$, combing with
$\frac{d^{2}s}{dt^{2}}|_{t=0}\geq0$, we have $\frac{d^{2}s}{dt^{2}}>0$ for $t>0$.
Taking into account to $\frac{ds}{dt}\mid_{t=0}<0$ and at the time $t_5$ when the delta shock wave disappears
$\frac{ds}{dt}|_{t=t_5}=u_-e^{-t_5}>0$, there exists a unique $t_9<t_5$ such that $\frac{ds}{dt}|_{t=t_9}=0$.
Moreover, $\frac{ds}{dt}<0$ for $0\leq t<t_9$ and $\frac{ds}{dt}>0$ for $t_9<t<t_5$. So the delta shock wave is convex for $0\leq t<t_5$, decreases along with $t$ such that the delta shock wave move fbackward for $0\leq t<t_9$, increases along with $t$ such that the delta shock wave move forward for $t_9<t<t_5$ and then disappears.
We draw Fig.4(b) to depict this situation,  where $t_9<t_2^\ast<t_5$.


\begin{multicols}{2}
\unitlength 1mm 
\linethickness{0.4pt}
\ifx\plotpoint\undefined\newsavebox{\plotpoint}\fi 
\begin{picture}(89.75,52.75)(10,0)
\put(89.25,7.25){\vector(1,0){.07}}
\put(13.75,7.25){\line(1,0){75.5}}
\put(89.75,3.75){$x$}
\put(87.75,20.75){\line(0,1){0}}
\put(52.5,51){\vector(0,1){.07}}
\put(52.5,7.25){\line(0,1){43.75}}
\qbezier(76.25,7.25)(9.38,21.63)(45,52.25)
\bezier{70}(33,32.75)(60,32.75)(77,32.75)
\qbezier(29,52.75)(33.25,14.38)(44.5,7.25)
\qbezier(33.25,28.25)(41.13,10.25)(52.5,7.25)
\qbezier(20,51.5)(23.75,13.5)(34.5,7.25)
\qbezier(41,14.75)(43.13,9.75)(47.75,7.25)
\qbezier(41.25,14.75)(47.38,9.75)(62,7.25)
\qbezier(51.25,51.75)(12.75,24.25)(84.25,7.25)
\put(32.5,48.75){Vac.}
\put(50,49.25){$t$}
\put(68,29){$t=t_2^\ast$}
\put(36,27.25){$t_5$}
\put(32.25,28){$\bullet$}
\put(32,31.5){$\bullet$}
\put(37,31.5){$\bullet$}
\put(51,4){0}
\put(35,12.5){$\delta S$}
\put(30,1){(a) $u_-+u_+<0$ for $u_+<u_-\leq0$}
\end{picture}


\unitlength 1mm 
\linethickness{0.4pt}
\ifx\plotpoint\undefined\newsavebox{\plotpoint}\fi 
\begin{picture}(90.75,51.5)(10,0)
\put(90.75,9.5){\vector(1,0){.07}}
\put(5,9.65){\line(1,0){85.75}}
\put(46.75,51.5){\vector(0,1){.07}}
\put(46.75,9.65){\line(0,1){42}}
\qbezier(17.25,9.65)(51,26.88)(59.75,46.75)
\qbezier(67,45.5)(29.75,23.38)(64.5,9.65)
\bezier{60}(48.25,23.25)(70,23.25)(87,23.25)
\qbezier(75.75,44.75)(33.5,22.75)(73.25,9.65)
\qbezier(51.5,33.5)(32.75,12.38)(47,9.65)
\qbezier(40.75,15)(35.25,9.75)(27.75,9.65)
\qbezier(41,15)(44.75,11.75)(51.5,9.65)
\qbezier(52.25,46)(49.38,37)(10,9.65)
\qbezier(42.25,20.25)(42.13,14.5)(55.5,9.65)
\qbezier(42,20.25)(37.75,15.5)(25.5,9.75)
\put(44,48.5){$t$}
\put(90.25,5.5){$x$}
\put(62,46){Vac.}
\put(53,32.5){$t_5$}
\put(49.5,31.5){$\bullet$}
\put(47,22){$\bullet$}
\put(53.5,22){$\bullet$}
\put(37,14){$t_9$}
\put(40,14){$\bullet$}
\put(45,29){$\bar{t}$}
\put(45.75,26.5){$\bullet$}
\put(81.25,20){$t=t_2^\ast$}
\put(47.25,6){0}
\put(20,12.5){$\delta S$}
\put(20,3){(b) $u_-+u_+<0$ for $u_+<0<u_-\leq1$}
\end{picture}
\put(-80,-6){\makebox(0,0)[cc]
{Fig.4 The delta shock wave solution to (1.1) and (1.2) when $f(x,t,u)=1$ }}
\put(-80,-10){\makebox(0,0)[cc]
{and $g(x,t,u)=-u$ for $u_-+u_+<0$ and $u_-\leq1$}}

\end{multicols}
\bigbreak

(ii)\hspace{0.15cm} \hangafter 1
\hangindent 3.5em
\noindent
If $u_{-}>1$, then $\frac{d^{2}s}{dt^{2}}|_{t=0}=\frac{1}{2}(1-u_-)<0$. Taking into account to
$\lim\limits{t\rightarrow+\infty}\frac{d^{2}s}{dt^{2}}=\frac{1}{2}>0$ and $\frac{d^{3}s}{dt^{3}}=\frac{1}{2}u_{-}e^{-t}>0$, we can conclude that there exists a unique $t_6$ such that $\frac{d^{2}s}{dt^{2}}|_{t=t_6}=0$, with $\frac{d^{2}s}{dt^{2}}<0$ for $0<t<t_6$ and $\frac{d^{2}s}{dt^{2}}>0$ for $t>t_6$.
Combing with $\frac{ds}{dt}|_{t=0}=\frac{1}{2}(u_-+u_+)<0$, we have $\frac{ds}{dt}|_{t=t_6}=\frac{1}{2}(u_-e^{-t_6}+u_++t_6)<\frac{ds}{dt}|_{t=0}<0$.
Similar to Case 1(ii)(b) in this section, we can prove $t_6<t_5$, which means that the delta shock wave disappears after the time $t_6$.
Taking into account to $\frac{ds}{dt}|_{t=t_5}=u_-e^{-t_5}>0$ and $\frac{d^{2}s}{dt^{2}}>0$ for $t_6<t<t_5$, there exists a unique $t_9$ such that $\frac{ds}{dt}|_{t=t_9}=0$ and $t_6<t_9<t_5$.
Moreover, $\frac{ds}{dt}<0$ for $0\leq t<t_9$ and $\frac{ds}{dt}>0$ for $t_9<t<t_5$.
So the delta shock wave is concave for $0< t<t_6$, convex for $t_6< t<t_5$, decreases along with $t$ such that the delta shock wave move fbackward for $0\leq t<t_9$, increases along with $t$ such that the delta shock wave move forward for $t_9<t<t_5$ and then disappears. We draw Fig.5 to depict this situation, where $t_9<t_2^\ast<t_5$.


\unitlength 1mm 
\linethickness{0.4pt}
\ifx\plotpoint\undefined\newsavebox{\plotpoint}\fi 
\begin{picture}(90.75,52.75)(-20,0)
\put(90.75,9.5){\vector(1,0){.07}}
\put(5,9.45){\line(1,0){85.75}}
\put(44,48.5){$t$}
\put(90.25,5.5){$x$}
\put(47.25,6){0}
\qbezier(12.75,9.45)(46.5,21.63)(68.25,45)
\put(47,52.75){\vector(0,1){.07}}
\put(47,10){\line(0,1){42.75}}
\qbezier(77.5,43.25)(27.88,25.75)(73.75,9.45)
\qbezier(84,41.75)(34.13,26.13)(80.75,9.45)
\bezier{50}(52,24.5)(70,24.5)(88,24.5)
\qbezier(57,34)(36,18.38)(42,15.25)
\qbezier(42,15.25)(46.5,12.88)(47,9.45)
\qbezier(40.5,17.25)(34.38,13)(22.75,9.45)
\qbezier(40.75,17.25)(49.5,10.75)(57.25,9.45)
\qbezier(60,45.25)(53,33.63)(6,9.45)
\put(70,43.5){Vac.}
\put(81,20.75){$t=t_2^\ast$}
\put(54.25,35){$t_5$}
\put(56,32.5){$\bullet$}
\put(37.75,17.75){$t_9$}
\put(39.75,16){$\bullet$}
\put(49.25,24){$\bar{t}$}
\put(45.95,24.5){$\bullet$}
\put(50.95,23.5){$\bullet$}
\put(57.5,23.5){$\bullet$}
\put(40,12.25){$t_6$}
\put(41,14){$\bullet$}
\put(60,1){\makebox(0,0)[cc]
{Fig.4 The delta shock wave solution to (1.1) and (1.2) when $f(x,t,u)=1$ }}
\put(60,-2){\makebox(0,0)[cc]
{and $g(x,t,u)=-u$ for $u_-+u_+<0$ and $u_->1$}}

\end{picture}


 \bigbreak
\noindent Case 3 \hspace{0.25cm} If ${u_-}+u_+=0$, taking into account ${u_-}>u_+$, it is easy to get $u_->0>u_+$.  According to (\ref{4.12}),
$\frac{d^{2}s}{dt^{2}}|_{t=0}=\frac{1}{2}(1-u_-)$, there are two subcases needed to be considered as follows.
\bigbreak

(i)\hspace{0.15cm} \hangafter 1
\hangindent 3.5em
\noindent
If $0<u_{-}\leq1$, then $\frac{d^{2}s}{dt^{2}}|_{t=0}=\frac{1}{2}(1-u_-)\geq0$. Since $\frac{d^{3}s}{dt^{3}}=\frac{1}{2}u_{-}e^{-t}>0$,
we have $\frac{d^{2}s}{dt^{2}}>0$ for $t>0$.
Taking into account $\frac{ds}{dt}\mid_{t=0}=0$, we have $\frac{ds}{dt}>0$ for $t>0$. But the delta shock wave disappear at $t_5$, so we have
 $\frac{ds}{dt}>0$ and $\frac{d^{2}s}{dt^{2}}>0$ for $0<t<t_5$. So the delta shock wave is convex for $0\leq t<t_5$, and increases along with $t$ such that the delta shock wave move forward for $0<t<t_5$ and then disappears. For the delta shock wave curve, one see Fig.6(b).

\bigbreak

(ii)\hspace{0.15cm} \hangafter 1
\hangindent 3.5em
\noindent
If $u_{-}>1$, then $\frac{d^{2}s}{dt^{2}}|_{t=0}=\frac{1}{2}(1-u_-)<0$. Since next the discussion is similar to case 2 (ii), one can also refer to Fig.7 for the delta shock wave curve. So we omit it.

\bigbreak
\begin{rem}\label{rem:4.1}  Similarly, we can consider other situations combined by $f(x,t,u)=\pm1$ and $g(x,t,u)=\pm u$.
Different from the situation for $f(x,t,u)\neq g(x,t,u)$, it is easy to prove that the delta shock wave always exists and never disappears
in the situation for $f(x,t,u)=g(x,t,u)=\pm1$ or $\pm u$.
\end{rem}

\section{Disccusions and Conclusions }
In this paper, we have finished the constructions of the delta shock wave solution to the Riemann problem (\ref{1.1}) and (\ref{1.2}) in all cases provided the
delta shock wave exists for the discontinuous source term $H(x-s(t))f(x,t,u)+H(s(t)-x)g(x,t,u)$ with special choices of $f(x,t,u)$ and $g(x,t,u)$. It is clear to see that .

It is noticed that the special choice of the discontinuous source term $H(x-s(t))f(x,t,u)+H(s(t)-x)g(x,t,u)$ is just for convenient in this paper. It is expected to
adopt the more
general discontinuous source term $H(x)f(x,t,u)+H(-x)g(x,t,u)$ in our later work.
Moreover, it is hoped that the method developed here can be used to study the pressureless Euler equations or the shallow water equations with singular source term.
Furthermore, we expect that this paper will give us some value insights into our future challenging study about the Cauchy problem for the pressureless Euler equations or the shallow water equations with source term.

\end{document}